\documentclass[11pt]{amsart}
\usepackage[margin=1in]{geometry}
\usepackage{amsmath,amssymb,amsthm}
\usepackage{thmtools}
\usepackage{xcolor}
\usepackage{tikz}
\usetikzlibrary{arrows.meta}
\usepackage[colorlinks=true,linkcolor=blue,citecolor=blue,urlcolor=blue]{hyperref}
\usepackage[capitalise,noabbrev,nameinlink]{cleveref}

\theoremstyle{plain}
\newtheorem{theorem}{Theorem}

\newtheorem{lemma}[theorem]{Lemma}
\newtheorem{proposition}[theorem]{Proposition}
\newtheorem{corollary}[theorem]{Corollary}
\theoremstyle{definition}
\newtheorem{definition}[theorem]{Definition}

\theoremstyle{remark}
\newtheorem{remark}[theorem]{Remark}

\crefname{theorem}{Theorem}{Theorems}
\crefname{maintheorem}{Theorem}{Theorems}
\crefname{lemma}{Lemma}{Lemmas}
\crefname{proposition}{Proposition}{Propositions}
\crefname{corollary}{Corollary}{Corollaries}
\crefname{definition}{Definition}{Definitions}
\crefname{example}{Example}{Examples}
\crefname{remark}{Remark}{Remarks}
\crefname{section}{Section}{Sections}

\newcommand{\R}{\mathbb{R}}
\newcommand{\Z}{\mathbb{Z}}
\newcommand{\VRle}[2]{\operatorname{VR}_{\le}(#1;#2)}
\newcommand{\VRlt}[2]{\operatorname{VR}_{<}(#1;#2)}
\newcommand{\Blt}[3]{B_{<}(#1,#2;#3)}
\newcommand{\Ble}[3]{B_{\le}(#1,#2;#3)}
\newcommand{\Tlt}[3]{T_{<}(#1,#2;#3)}
\newcommand{\Tle}[3]{T_{\le}(#1,#2;#3)}

\title[Vietoris--Rips Contractibility in Normed Spaces]{Vietoris--Rips
Contractibility for Dense Subsets of Finite-Dimensional Normed Spaces}

\author{Qingsong Wang}
\address{Hal{\i}c{\i}o\u{g}lu Data Science Institute,
University of California San Diego, La Jolla, CA 92093, USA}
\email{qswang92@gmail.com}

\author{Ling Zhou}
\address{Department of Mathematics,
University of Oklahoma, Norman, OK 73019, USA}
\email{zhouling0903@gmail.com}

\date{August 15, 2026}

\begin{document}

\begin{abstract}
Let \(Y\) be a \(\delta\)-dense subset of a finite-dimensional normed space
\(X\), so every point of \(X\) is at distance at most \(\delta\) from some
point of \(Y\). We prove that the open Vietoris--Rips complex
\(\VRlt{Y}{r}\) is contractible whenever
\(r>\delta/(1-J_{\mathrm{fin}}(X))\). Here \(J_{\mathrm{fin}}(X)\) is the
normalized finite Jung constant, namely the supremum of circumradius divided
by diameter over finite subsets of positive diameter.
To prove this, we use the canonical isometric embedding of \(X\) into a
hyperconvex normed space. There \(\VRlt{Y}{r}\) is the nerve of an equal-radius
open-ball cover. At the stated scales, the finite Jung estimate allows the
cover to be restricted to a convex tube about \(X\) without changing its nerve.
The \(\delta\)-density condition then makes the restricted balls cover the
entire tube. Consequently, \(\VRlt{Y}{r}\) is homotopy equivalent to the tube
and hence contractible.

Applying this result to lattices gives explicit contractibility bounds. In
\((\Z^n,\ell_1)\), gaps in the distance spectrum further sharpen the resulting
bound: the open Vietoris--Rips complex \(\VRlt{\Z^n,\ell_1}{r}\) is contractible above \(\max\{n,n(n-1)/2\}\), and the
closed Vietoris--Rips complex \(\VRle{\Z^n,\ell_1}{r}\) is contractible at or above the same value
due to its distance spectrum being discrete. This asymptotically
halves the leading term of Zaremsky's previous bound. 
Our method also gives explicit bounds for \((\Z^n,\ell_p)\) for all $1\leq p \leq \infty$.
In particular, it recovers the optimal bound \(r=1\) when \(p=\infty\). Meanwhile, for fixed \(n\), the bounds converge when letting \(p\to \infty\). For a general locally finite \(\delta\)-dense subset
\(Y\subset X\) whose distance spectrum may not be discrete, a separate
closed-cover argument proves that \(\VRle{Y}{r}\) is contractible whenever
\(r\ge\delta/(1-J_{\mathrm{fin}}(X))\).
\end{abstract}

\maketitle

\section{Introduction}
The Vietoris--Rips construction associates to each metric space a family of
simplicial complexes indexed by a scale parameter. It is a basic tool in
metric geometry and applied topology~\cite{gromov,chazal2009}. 
For a metric
space \((M,d_M)\) and \(r>0\), let
\(\VRle{M,d_M}{r}\) be the simplicial complex whose simplices are the finite
subsets \(\sigma\subset M\) satisfying
\(\operatorname{diam}(\sigma)\le r\), and define
\(\VRlt{M,d_M}{r}\) by replacing \(\le r\) with \(<r\). 
These are the closed
and open Vietoris--Rips complexes, respectively. When the metric is clear, we
write \(\VRle{M}{r}\) and \(\VRlt{M}{r}\).

If \(M\) has finite diameter, then \(\VRle{M}{r}\) is a simplex for
\(r\ge\operatorname{diam}(M)\) and hence is contractible. For an unbounded
space, no finite scale joins every pair of vertices, so contractibility must
arise from additional geometry. A finitely generated group is
\emph{word-hyperbolic} if its Cayley graph with respect to one (equivalently,
every) finite generating set is Gromov hyperbolic. For such a group \(G\), the
closed Vietoris--Rips complex \(\VRle{G}{r}\) is contractible for all sufficiently
large \(r\)~\cite[Section~2.2]{gromov}. For the standard generators of
\(\Z^n\), the word metric is the restriction of the \(\ell_1\)-metric. Since
\((\Z^n,\ell_1)\) is not Gromov hyperbolic for \(n\ge2\), the
hyperbolic-group result does not apply.

{
More broadly, it is interesting to understand how far
large-scale Rips contractibility extends beyond hyperbolic groups. The lattice
\((\Z^n,\ell_1)\) is a basic nonhyperbolic test case: Zaremsky asked whether
\(\VRle{\Z^n,\ell_1}{r}\) is contractible at all sufficiently large scales
~\cite[Section~6.1]{zaremskybb} and conjectured that
\(\VRle{\Z^n,\ell_1}{r}\) is contractible for all
\(n\ge1\) and \(r\ge n\)~\cite{zaremsky}.
The lower bound \(r\ge n\) is necessary, since
\(\VRle{\Z^n,\ell_1}{r}\) is not contractible for \(r<n\);
see \cite[Section~6]{virk}.
Virk proved the optimal 
contractibility for \(r\ge n\) when \(n\le3\) and for
\(r\ge n^2(2n-1)\) in every dimension~\cite{virk}.
Later, Zaremsky obtained the bound \(n^2+n-1\), which was previously the best
all-dimensional estimate~\cite{zaremsky}.
Meanwhile, McCarty's dissertation gives the related local-domination
bound \(r\ge n(n+1)\)~\cite[Theorem~3.1]{mccarty}. Gupta, Sarkar, and Shukla
proved Zaremsky's conjecture for \(n\le5\) and established contractibility
for \(\Z^6\) when \(r\ge10\)
\cite{guptasarkarshukla}.
\par}

{In this paper, we improve Zaremsky's all-dimensional
contractibility bound for \((\Z^n,\ell_1)\) from \(r\ge n^2+n-1\) to
\(r\ge\max\{n,n(n-1)/2\}\), asymptotically halving its leading term.
Our general theorem gives contractibility bounds for \(\delta\)-dense subsets of
finite-dimensional normed spaces. For lattices, explicit density calculations
make these bounds concrete, and quantitative gaps in the lattice distance spectrum yield the sharper \(\ell_1\) estimate.
\par}

{
Let \(X\) be a finite-dimensional normed space and let \(Y\subset X\) be
\(\delta\)-dense, meaning that every point of \(X\) is at distance at most
\(\delta\) from some point of \(Y\).
Let \(J_{\mathrm{fin}}(X)\), called the \emph{(normalized) finite Jung constant}, denote the
supremum of the ratios of circumradius to diameter
over all finite subsets of \(X\) with positive diameter.
If \(\dim X=n\ge1\), every two-point subset of \(X\) has circumradius equal to
half its diameter, while Bohnenblust's inequality~\cite[Theorem~6]{bohnenblust}
gives the upper bound
\[
\frac12\le J_{\mathrm{fin}}(X)\le\frac{n}{n+1}<1.
\]
For \(X=\{0\}\), one has \(J_{\mathrm{fin}}(X)=0\).
Hence \(1-J_{\mathrm{fin}}(X)>0\) for every finite-dimensional normed space
\(X\). Our main open-convention theorem is the following.
}

\begin{restatable}[Contractibility of open Vietoris--Rips complexes]{maintheorem}{openjungcontractibility}
\label{thm:jung-open-rips}
Let \(X\) be a finite-dimensional normed space, and let \(Y\subset X\) be
\(\delta\)-dense.
Then \(\VRlt{Y}{r}\) is contractible whenever
\[
r>\frac{\delta}{1-J_{\mathrm{fin}}(X)}.
\]
\end{restatable}

{\paragraph{\textbf{Hyperconvex embeddings.}}

Our proof uses a result of Lim, Mémoli, and Okutan relating open
Vietoris--Rips complexes to tubular neighborhoods in \emph{hyperconvex}
metric spaces~\cite[Propositions~2.25 and~2.27]{limmemoliokutan}.
Hyperconvex metric spaces were introduced by Aronszajn and Panitchpakdi,
who proved that they are precisely the \emph{injective} metric spaces
\cite{aronszajnpanitchpakdi}. 
Fundamental examples include spaces of bounded real-valued functions on an
index set \(I\), equipped with the supremum norm and denoted by
\(\ell_\infty(I)\).
Every metric space embeds isometrically into one of these
spaces. 
Isbell proved that among the injective extensions of a metric space
there is a minimal one, unique up to isometry over the original space, called
its \emph{injective envelope}~\cite{isbell}. 
Dress later gave an equivalent
construction under the name \emph{tight span}~\cite{dress}.
Subsequent work has connected injective hulls with
{finite metric geometry and split
decompositions~\cite{bandeltdress}, phylogenetic
combinatorics~\cite{dresshuberkoolenmoultonspillner},}
discrete metric spaces and groups~\cite{herrlich,lang},
{and absolute \(1\)-Lipschitz retracts in
\(\ell_\infty(I)\)~\cite{descombesspavon}.}

The result of Lim, Mémoli, and Okutan used in our proof can be stated as
follows: if \(Y\) embeds isometrically into a hyperconvex
(equivalently, injective) metric space \(H\), then \(\VRlt{Y}{r}\) is
homotopy equivalent to the open \(r/2\)-neighborhood of the embedded copy
of \(Y\) in \(H\)
\cite[Propositions~2.25 and~2.27]{limmemoliokutan}. This result extends
earlier nerve identifications for finite metric spaces in their canonical
\(\ell_\infty\)-embeddings: see~\cite[Lemma~4]{ghristmuhammad} for closed
cubes and \cite[Lemma~2.9]{chazal2009} for open balls.

In our paper, for \(Y\subset X\), we consider the following (linear) hyperconvex embedding of normed vector space:
\[
\iota:X\hookrightarrow H:=\ell_\infty(B_{X^*}),
\]
where \(B_{X^*}\) is the closed unit ball of the dual space \(X^*\). We identify \(X\) with
its image and regard \(Y\subset X\subset H\). 

The comparison identifies \(\VRlt{Y}{r}\), up to homotopy, with the open
\(r/2\)-neighborhood of \(Y\) in \(H\).
At the scales in \cref{thm:jung-open-rips}, the \(\delta\)-density of
\(Y\) allows us to choose an open tube around \(X\) that is covered by the
open \(r/2\)-balls centered at \(Y\).
Since \(X\) is a linear subspace of \(H\), this tube is convex
and therefore contractible.
It remains to compare the original ball cover with
its restriction to the tube. 
In \cref{prop:finite-jung-tube-restriction}, we use finite Jung estimates
to show that the two covers have the same nerve. The estimates used in our
applications come from classical results of Bohnenblust for arbitrary
finite-dimensional normed spaces~\cite{bohnenblust} and of Ivanov and
Pichugov for \(\ell_p^n\)~\cite{ivanovpichugov}.
Since the union of the restricted
balls is the tube, {the equality of the
two nerves completes the argument: \(\VRlt{Y}{r}\) is homotopy equivalent
to a convex tube and is therefore contractible.}

The proof mechanism is summarized in
\cref{fig:ball-tube-schematic}.

\begin{figure}[htbp]
\centering
\begin{tikzpicture}[font=\small]
\definecolor{YBlue}{HTML}{0072B2}
\definecolor{XGold}{HTML}{E69F00}
\definecolor{RestrictionGreen}{HTML}{009E73}
\tikzset{
  ambient/.style={draw=black!45, rounded corners=2pt, fill=black!2},
  ball/.style={draw=YBlue, densely dashed, fill=YBlue,
    fill opacity=0.24},
  ballboundary/.style={draw=YBlue, densely dashed, line width=0.8pt},
  restrictedball/.style={draw=RestrictionGreen, densely dashed,
    fill=RestrictionGreen!24},
  restrictedboundary/.style={draw=RestrictionGreen, densely dashed,
    line width=0.5pt},
  tube/.style={draw=XGold, dashed, fill=XGold,
    fill opacity=0.18},
  xspace/.style={draw=XGold, very thick},
  ypoint/.style={circle, fill=YBlue, inner sep=1.5pt}
}

\begin{scope}
  \draw[ambient] (0,0) rectangle (4.2,3);
  \node[anchor=north west] at (0.12,2.88) {\(H\)};
  \draw[ball] (0.28,0.37) rectangle (1.32,1.41);
  \draw[ball] (1.13,0.77) rectangle (2.17,1.81);
  \draw[ball] (1.98,1.17) rectangle (3.02,2.21);
  \draw[ball] (2.83,1.57) rectangle (3.87,2.61);
  \draw[xspace] (0.3,0.65) -- (3.9,2.35)
    node[pos=0.97, right, text=XGold] {\(X\)};
  \foreach \x/\y in {0.8/0.89,1.65/1.29,2.5/1.69,3.35/2.09}
    \node[ypoint] at (\x,\y) {};
  \node[YBlue] at (2.1,0.16) {ball cover centered at \(Y\)};
  \node[font=\bfseries] at (2.1,3.3)
    {\textcolor{YBlue}{(a)} Ambient ball cover};
\end{scope}

\begin{scope}[xshift=5.4cm]
  \draw[ambient] (0,0) rectangle (4.2,3);
  \node[anchor=north west] at (0.12,2.88) {\(H\)};
  \path[tube] (0.2,0.22) -- (4,2.02) -- (4,2.78) -- (0.2,0.98) -- cycle;
  \begin{scope}[even odd rule]
    \clip
      (-0.01,-0.01) rectangle (4.21,3.01)
      (0.2,0.22) -- (4,2.02) -- (4,2.78) -- (0.2,0.98) -- cycle;
    \draw[ballboundary] (0.28,0.37) rectangle (1.32,1.41);
    \draw[ballboundary] (1.13,0.77) rectangle (2.17,1.81);
    \draw[ballboundary] (1.98,1.17) rectangle (3.02,2.21);
    \draw[ballboundary] (2.83,1.57) rectangle (3.87,2.61);
  \end{scope}
  \begin{scope}
    \clip (0.2,0.22) -- (4,2.02) -- (4,2.78) -- (0.2,0.98) -- cycle;
    \draw[restrictedball] (0.28,0.37) rectangle (1.32,1.41);
    \draw[restrictedball] (1.13,0.77) rectangle (2.17,1.81);
    \draw[restrictedball] (1.98,1.17) rectangle (3.02,2.21);
    \draw[restrictedball] (2.83,1.57) rectangle (3.87,2.61);
    \draw[restrictedboundary] (0.28,0.37) rectangle (1.32,1.41);
    \draw[restrictedboundary] (1.13,0.77) rectangle (2.17,1.81);
    \draw[restrictedboundary] (1.98,1.17) rectangle (3.02,2.21);
    \draw[restrictedboundary] (2.83,1.57) rectangle (3.87,2.61);
  \end{scope}
  \draw[XGold, dashed]
    (0.2,0.22) -- (4,2.02) -- (4,2.78) -- (0.2,0.98) -- cycle;
  \draw[xspace] (0.3,0.65) -- (3.9,2.35)
    node[pos=0.97, right, text=XGold] {};
  \foreach \x/\y in {0.8/0.89,1.65/1.29,2.5/1.69,3.35/2.09}
    \node[ypoint] at (\x,\y) {};
  \node[inner sep=2pt, text=XGold]
    at (3,2.80) {convex tube \(T\)};
  \node[text=RestrictionGreen] at (2.1,0.16) {restricted ball cover};
  \node[font=\bfseries] at (2.1,3.3)
    {\textcolor{RestrictionGreen}{(b)} Restriction to the tube};
\end{scope}

\begin{scope}[xshift=10.8cm]
  \draw[ambient] (0,0) rectangle (4.2,3);
  \node[anchor=north west] at (0.12,2.88) {\(H\)};
  \path[draw=XGold, dashed, fill=RestrictionGreen!24]
    (0.2,0.22) -- (4,2.02) -- (4,2.78) -- (0.2,0.98) -- cycle;
  \draw[xspace] (0.3,0.65) -- (3.9,2.35)
    node[pos=0.97, right, text=XGold] {};
  \foreach \x/\y in {0.8/0.89,1.65/1.29,2.5/1.69,3.35/2.09}
    \node[ypoint] at (\x,\y) {};
  \node at (2.1,0.16)
    {\textcolor{RestrictionGreen}{restricted union}
      \(=\textcolor{XGold}{T}\)};
  \node[font=\bfseries] at (2.1,3.3)
    {\textcolor{XGold}{(c)} Density fills the tube};
\end{scope}

\node[font=\normalsize] at (7.5,-0.62)
  {\(\VRlt{\textcolor{YBlue}{Y}}{2\rho}
    \overset{\textcolor{YBlue}{\text{(a)}}}{\simeq}
    \textcolor{YBlue}{\text{ball union}}
    \overset{\textcolor{RestrictionGreen}{\text{(b)}}}{\simeq}
    \textcolor{RestrictionGreen}{\text{restricted union}}
    \overset{\textcolor{XGold}{\text{(c)}}}{=}
    \textcolor{XGold}{T}\simeq *\)};
\end{tikzpicture}
\caption{Geometric outline of the contractibility argument. \textup{(a)} In
the hyperconvex ambient space \(H\), the open \(\rho\)-ball cover centered at
\(Y\) has nerve \(\VRlt{Y}{2\rho}\). \textup{(b)} We use the finite Jung bound to show that restricting this cover to a tube \(T\) about \(X\) preserves its nerve.
\textup{(c)} The density of \(Y\) allows us to identify the union of the
restricted cover with the convex tube \(T\). Hence
\(\VRlt{Y}{2\rho}\simeq T\simeq *\).}
\label{fig:ball-tube-schematic}
\end{figure}

}

\medskip

We now consider the contractibility of Vietoris--Rips complexes of lattices in normed spaces.
Specifically, we focus on full-rank \emph{affine lattices} of the form
\(\Lambda=a+\sum_{i=1}^n\Z v_i\), where \(a\in X\) and
\(v_1,\ldots,v_n\) is a basis of \(X\).
When \(a=0\), we call
\(\Lambda\) a full-rank \emph{lattice}.
Compactness of a fundamental parallelepiped implies that \(\Lambda\) is
\(\delta\)-dense for some finite \(\delta\).
We can then apply \Cref{thm:jung-open-rips} to obtain contractibility for the (open) Vietoris-Rips $\VRlt{\Lambda}{s}$. 
Additionally, we can extend the result to closed Vietoris--Rips complexes. 
This is because the distance spectrum, namely the set of pairwise distances in \(\Lambda\), is
locally finite. Hence, for every \(r>0\), there exists \(s>r\) such that no
pairwise lattice distance lies in \((r,s]\), and therefore
\[
\VRle{\Lambda}{r}=\VRlt{\Lambda}{s}.
\]
Consequently, a bound \(R\) for all open scales \(s>R\) gives the same bound
for all closed scales \(r\ge R\); see
\cref{prop:open-closed-threshold-equivalence}.
In summary, we have the following result. 
\begin{restatable}[Lattices in normed spaces]{maintheorem}{normedlatticecontractibility}
\label{thm:normed-lattice}
{
Let \(X\) be a finite-dimensional normed space, and let \(\Lambda\subset X\)
be a full-rank affine lattice that is \(\delta\)-dense. Then
\(\VRlt{\Lambda}{r}\) is contractible whenever
\[
r>\frac{\delta}{1-J_{\mathrm{fin}}(X)},
\]
and \(\VRle{\Lambda}{r}\) is contractible whenever
\[
r\ge\frac{\delta}{1-J_{\mathrm{fin}}(X)}.
\]
}
\end{restatable}

Both \(\delta\) and \(J_{\mathrm{fin}}(X)\) are computed with respect to the
chosen norm on \(X\). Rescaling the norm rescales
\(\delta\) and the Vietoris--Rips parameter but leaves
\(J_{\mathrm{fin}}(X)\) unchanged, so comparisons between norms require a
fixed normalization. We use the standard normalization of the \(\ell_p\)
norms, in which the coordinate vectors have norm one.

Let \(\ell_p^n\) denote \(\R^n\) equipped with the \(\ell_p\)-norm.
To apply \cref{thm:normed-lattice} to
\((\Z^n,\ell_p)\), we need a value of \(\delta\) for which
\(\Z^n\) is \(\delta\)-dense in \(\ell_p^n\), together with an upper bound for \(J_{\mathrm{fin}}(\ell_p^n)\). Coordinatewise rounding shows that
\(\Z^n\) is \(n^{1/p}/2\)-dense. The following result combines this density estimate with the required Jung estimate and \cref{thm:normed-lattice}.

\begin{restatable}[The integer lattice with the \(\ell_p\) metric]{maintheorem}{lplatticecontractibility}
\label{thm:lp-lattice}
\label{cor:lp-lattice}
{
Let \(n\ge1\) and \(1\le p<\infty\). Define \(p'\in(1,\infty]\) by
\(1/p+1/p'=1\), where \(1/\infty:=0\), and set
\(s_p:=\max\{p,p'\}\). Define
\begin{equation}
\label{eq:j_{p,n}}
j_{p,n}:=\min\left\{
2^{-1/s_p}\left(\frac{n}{n+1}\right)^{1/p},
\frac{n^{1/p}}2
\right\}.
\end{equation}
Set
\begin{equation}
\label{eq:r_{p,n}}
r_{p,n}:=\frac{n^{1/p}}{2(1-j_{p,n})}.
\end{equation}
Then \(\VRlt{\Z^n,\ell_p}{r}\) is contractible for
\(r>r_{p,n}\), and \(\VRle{\Z^n,\ell_p}{r}\) is contractible for
\(r\ge r_{p,n}\).
For \(p=\infty\), the open complex is contractible if and only if \(r>1\),
and the closed complex is contractible if and only if \(r\ge1\).
}
\end{restatable}

For
fixed \(n\), the scales $r_{p,n}$ converge as \(p\to\infty\) to the exact
\(\ell_\infty\) threshold \(r_{\infty,n}=1\). When \(p\) is a positive integer, the
\(p\)-th powers of \(\ell_p\)-distances between lattice points are integers.
The resulting distance gaps strictly lower the sufficient scale for
\(n\ge2\); see
\cref{cor:integer-p-refinement}. In dimensions one and two, the lattice
structure determines the exact open and closed thresholds; see
\cref{cor:lp2-exact-threshold}.

{
For \(p=1\), one has \(j_{1,n}=n/(n+1)\) and
\(r_{1,n}=n(n+1)/2\).
The integer-valued \(\ell_1\)-distance spectrum allows us to improve this
contractibility threshold, as follows.
}

\begin{restatable}[The integer lattice with the \(\ell_1\) metric]{corollary}{integerlatticecontractibility}
\label{thm:integer-lattice}
\label{cor:lattice}
{
For every \(n\ge1\), set
\[
K_n:=\max\left\{n,\frac{n(n-1)}2\right\}.
\]
Then \(\VRlt{\Z^n,\ell_1}{r}\) is contractible for \(r>K_n\), and
\(\VRle{\Z^n,\ell_1}{r}\) is contractible for \(r\ge K_n\).
}
\end{restatable}

{
Together with Virk's lower bound~\cite{virk}, the equality \(K_n=n\) for
\(n\le3\) shows that \cref{cor:lattice} is sharp in dimensions one, two, and
three. For \(n\ge4\), one has
\(K_n=n(n-1)/2\); this asymptotically halves the leading term of Zaremsky's
bound \(n^2+n-1\)~\cite{zaremsky}.
Our bound comes from a uniform-tube argument for general
dense subsets. We suspect that a tube adapted to the lattice geometry may
reach the conjectural threshold \(r=n\).
}

\medskip

{
The passage from \cref{thm:jung-open-rips} to \cref{thm:normed-lattice} uses
the locally finite distance spectrum of the lattice. A locally finite
\(\delta\)-dense subset need not have a locally finite distance spectrum: its
pairwise distances can accumulate in a bounded interval. Thus
\cref{prop:open-closed-threshold-equivalence} need not apply. We instead use
Nag\'orko's closed-cover nerve theorem to prove the closed
result for every locally finite \(\delta\)-dense subset of a finite-dimensional
normed space. Local finiteness and convexity of the ball intersections provide
the cover hypotheses~\cite[Theorem~3.3]{nagorko}, while finite dimensionality
gives the attainment needed at the endpoint.
{
A subset \(Y\subset X\) is \emph{locally finite} if every compact subset of
\(X\) meets \(Y\) in finitely many points.
}
}

\begin{restatable}[Contractibility of closed Vietoris--Rips complexes]{maintheorem}{closedellonecontractibility}
\label{thm:main-rips}
Let \(X\) be a finite-dimensional normed space, and let \(Y\subset X\) be
locally finite and \(\delta\)-dense.
Then \(\VRle{Y}{r}\) is contractible whenever
\[
r\ge\frac{\delta}{1-J_{\mathrm{fin}}(X)}.
\]
\end{restatable}

{
In dimension two, the isometry
\(\ell_1^2\cong\ell_\infty^2\) gives
\(J_{\mathrm{fin}}(\ell_1^2)=1/2\). Thus the general open and closed theorems
specialize to the bounds \(r>2\delta\) and \(r\ge2\delta\), respectively.
See \cref{cor:open-l12-sharp,cor:closed-planar-lp}.}

\subsection*{Organization.}
{
\Cref{sec:preliminaries} fixes the metric notation and finite Jung constant.
\Cref{sec:hyperconvexity-relative-density} proves
\cref{thm:jung-open-rips}. It also derives explicit open
bounds for \(\delta\)-dense subsets of \(\ell_p^n\).
\Cref{sec:lp-lattice-consequences}
computes the density of coordinate lattices, translates the open bounds to
closed endpoints, refines them using gaps in the realized distance spectrum,
and determines
the exact thresholds in dimensions one and two.
\Cref{sec:closed-convention} proves the closed theorem for locally finite
\(\delta\)-dense subsets and gives its planar specializations.
}

\subsection*{Acknowledgments}
We thank Facundo M\'emoli for suggesting this question to us and for
his encouragement. Q.W. also thanks Matthew Zaremsky for comments and
suggestions on an earlier version of the paper.

\section{Vietoris--Rips Complexes and Circumradii}
\label{sec:preliminaries}

All normed spaces in this paper are real and are equipped with the metric
induced by the norm; thus, for a normed space \(X\),
\(
d_X(x,x')=\|x-x'\|.
\)
We write \((M,d_M)\) for a generic
metric space; metrics carry their ambient-space subscripts, as in \(d_X\) and
\(d_M\).
When the same underlying set carries different metrics, the metric is included
in the Vietoris--Rips notation; for example,
\(\VRle{\Z^n,\ell_p}{r}\) denotes the closed Vietoris--Rips complex of
\(\Z^n\) with the metric induced by the \(\ell_p\)-norm. When the metric on
\(M\) has already been fixed, we abbreviate the Vietoris-Rips notation \(\VRle{M,d_M}{r}\) to \(\VRle{M}{r}\),
and likewise for the open convention. Diameters are taken with respect to the
specified metric, with a subscript added when more than one metric is present.
For \(1\le p<\infty\), the H\"older-conjugate exponent
\(p'\in(1,\infty]\) is defined by
\(1/p+1/p'=1\), where \(1/\infty:=0\).

\begin{definition}[Balls and tubes]
Let \((M,d_M)\) be a metric space and
\(\emptyset\ne A\subset M\). For any \(z\in M\), define {the \emph{distance from \(z\) to \(A\)} by}
\[
\operatorname{dist}_M(z,A):=\inf_{x\in A}d_M(z,x).
\]

Let \(a\in M\) and \(\rho\ge0\).
In the notation below, the center point \(a\) or set \(A\) comes first, the ambient
metric space \(M\) second, and the radius \(\rho\) follows the semicolon. The
\emph{open ball}, \emph{closed ball}, \emph{open tube}, and \emph{closed tube}
are defined by
\[
\begin{aligned}
\Blt{a}{M}{\rho}
&:=\{z\in M:d_M(z,a)<\rho\},
&
\Ble{a}{M}{\rho}
&:=\{z\in M:d_M(z,a)\le\rho\},\\
\Tlt{A}{M}{\rho}
&:=\{z\in M:\operatorname{dist}_M(z,A)<\rho\},
&
\Tle{A}{M}{\rho}
&:=\{z\in M:\operatorname{dist}_M(z,A)\le\rho\}.
\end{aligned}
\]
The subscripts \(<\) and \(\le\) indicate strict and weak inequalities, respectively;
in particular, \(B_{\le}\) does not denote the closure of \(B_{<}\).

We write
ball unions explicitly. One always has
\[
\bigcup_{a\in A}\Blt{a}{M}{\rho}=\Tlt{A}{M}{\rho},
\qquad
\bigcup_{a\in A}\Ble{a}{M}{\rho}\subseteq\Tle{A}{M}{\rho}.
\]
The second inclusion is an equality whenever the distance to \(A\) is attained
at every point of \(\Tle{A}{M}{\rho}\).
\end{definition}

\begin{definition}[\(\delta\)-dense subsets]
Let \((M,d_M)\) be a metric space. For \(A\subset M\) and \(\delta\ge0\), we say
that \(A\) is \emph{\(\delta\)-dense} in \(M\) if
\[
M=\bigcup_{a\in A}\Ble{a}{M}{\delta}.
\]
Equivalently, every \(x\in M\) lies in \(\Ble{a}{M}{\delta}\) for some
\(a\in A\).
\end{definition}

\begin{definition}[Diameter, circumradius and normalized finite Jung constant]
Let \((M,d_M)\) be a metric space. For every nonempty subset
\(A\subset M\), define its \emph{diameter} by
\[
\operatorname{diam}_M(A)
:=
\sup\{d_M(x,x'):x,x'\in A\}.
\]
For a finite nonempty set \(W\subset M\), define its \emph{circumradius} in
\(M\) by
\[
R_M(W)
:=
\inf_{c\in M}\max_{w\in W}d_M(c,w).
\]

The \emph{(normalized) finite Jung constant} of \(M\) is
\[
J_{\mathrm{fin}}(M)
:=
\sup_{\substack{W\subset M\text{ finite}\\
\operatorname{diam}_M(W)>0}}
\frac{R_M(W)}{\operatorname{diam}_M(W)}.
\]
Only finite subsets enter the definition because every Vietoris--Rips
simplex has finitely many vertices. If \(M\) has at most one point, we
set \(J_{\mathrm{fin}}(M)=0\). 
\end{definition}

Some Banach-space sources use the doubled normalization \(2R_X(W)/\operatorname{diam}(W)\). See \cite{castillopapini} for background and normalization conventions.

\begin{definition}[Vietoris--Rips complexes]
Let \((M,d_M)\) be a metric space and let \(r>0\). The \emph{closed
Vietoris--Rips complex} \(\VRle{M}{r}\) is the abstract simplicial complex with
vertex set \(M\) in which a nonempty finite subset \(\sigma\subset M\) spans a
simplex if and only if
\[
\operatorname{diam}_M(\sigma)\le r .
\]
The \emph{open Vietoris--Rips complex} \(\VRlt{M}{r}\) is defined by replacing
\(\le r\) with \(<r\). In both cases the empty simplex is included by
convention.
\end{definition}

When discussing homotopy type, we do not distinguish notationally
between an abstract simplicial complex and its geometric realization.

\section{Open Vietoris--Rips Complexes of
\texorpdfstring{\(\delta\)-Dense}{delta-Dense} Sets}
\label{sec:hyperconvexity-relative-density}

The canonical embedding
\(
X\hookrightarrow H:=\ell_\infty(B_{X^*})
\)
provides the ambient space for the proof of
\cref{thm:jung-open-rips}; see
\cref{subsec:hyperconvex-ambient}. In this hyperconvex space, the open
Vietoris--Rips complex is represented by an equal-radius open-ball union;
see \cref{subsec:open-ball-unions}. The finite Jung bound then allows this
ball cover to be restricted to an open tube about the embedded copy of \(X\)
without changing its nerve; see
\cref{subsec:open-tube-restriction}. At the scales specified in the theorem,
the \(\delta\)-density of \(Y\) makes the restricted balls cover the entire
tube, which is convex and hence contractible; see
\cref{subsec:open-tube-density}. Finally, the Jung estimates in
\cref{subsec:lp-open-estimates} give explicit bounds for
\(\delta\)-dense subsets of \(\ell_p^n\).

\subsection{Hyperconvex ambient spaces}
\label{subsec:hyperconvex-ambient}

\begin{definition}[Hyperconvexity]
A metric space \(H\) is \emph{hyperconvex} if every family of closed balls
\(\{\Ble{p_i}{H}{r_i}\}_{i\in I}\) satisfying
\[
d_H(p_i,p_j)\le r_i+r_j
\qquad(i,j\in I)
\]
has nonempty total intersection.

Hyperconvex metric spaces are precisely the
\emph{injective} metric spaces; this is the theorem of Aronszajn and
Panitchpakdi~\cite{aronszajnpanitchpakdi}, also stated as
\cite[Proposition~2.18]{limmemoliokutan}. 
We refer to these spaces as \emph{hyperconvex} and use only their defining
ball-intersection property.
\end{definition}

For every set \(I\), the Banach space \(\ell_\infty(I)\) of bounded functions
\(f:I\to\R\), equipped with the supremum norm, is injective
\cite[Example~2.22]{limmemoliokutan}, and hence hyperconvex. Thus
\(\ell_\infty\)-spaces turn pairwise compatibility of balls into a common
intersection. The Hahn--Banach theorem provides an isometric embedding of every
normed space into such an ambient space.

\begin{lemma}[Canonical \(\ell_\infty\)-embedding]
\label{lem:canonical-embedding}
Let \(X\) be a normed space, and let \(B_{X^*}\) be the closed unit ball of
the dual space $X^*$. The map
\[
\iota:X\longrightarrow \ell_\infty(B_{X^*}),
\qquad
\iota(x)(\varphi)=\varphi(x),
\]
is a linear isometric embedding.
\end{lemma}

\begin{proof}
The map \(\iota\) is linear. For every \(x\in X\) and
\(\varphi\in B_{X^*}\),
\[
|\iota(x)(\varphi)|=|\varphi(x)|
\le \|\varphi\|\,\|x\|\le\|x\|.
\]
Thus \(\iota(x)\in\ell_\infty(B_{X^*})\) and
\(\|\iota(x)\|_\infty\le\|x\|\).

If \(x\neq0\), the Hahn--Banach theorem provides
\(\varphi\in X^*\) such that \(\|\varphi\|=1\) and
\(\varphi(x)=\|x\|\). Hence
\[
\|\iota(x)\|_\infty
=\sup_{\psi\in B_{X^*}}|\psi(x)|
\ge|\varphi(x)|=\|x\|.
\]
The case \(x=0\) is immediate, so \(\iota\) is isometric.
\end{proof}

\subsection{Open ball unions in hyperconvex spaces}
\label{subsec:open-ball-unions}

An isometric embedding into a hyperconvex space realizes an open
Vietoris--Rips complex as the homotopy type of an open ball union.
Lim, M\'emoli, and Okutan identify the open Vietoris--Rips complex with
the nerve of the equal-radius open-ball cover, namely the corresponding
ambient \v{C}ech complex, and apply the open nerve theorem
\cite[Propositions~2.25 and~2.27]{limmemoliokutan}.
For finite metric spaces in
their canonical \(\ell_\infty\)-embeddings, the open-ball statement appears in
\cite[Lemma~2.9]{chazal2009}; \cite[Lemma~VII]{ghristmuhammad} gives the
corresponding closed-cube identity.

\begin{proposition}[Open Vietoris--Rips complexes as ball unions]
\label{prop:open-thickening}
\label{prop:open-ball-nerve}
Let \((M,d_M)\) be a metric space, let \(H\) be a hyperconvex
metric space, and let \(\eta:M\to H\) be an isometric embedding. For
\(\rho>0\), form the union of the open \(\rho\)-balls centered at
\(\eta(M)\).
Then
\[
\VRlt{M}{2\rho}
\simeq \bigcup_{m\in M}\Blt{\eta(m)}{H}{\rho}
=\Tlt{\eta(M)}{H}{\rho}.
\]
\end{proposition}

\subsection{Nerve-preserving restriction to a tube}
\label{subsec:open-tube-restriction}

The preceding proposition represents \(\VRlt{Y}{2\rho}\) by the union of the
open \(\rho\)-balls centered at \(Y\) in the hyperconvex ambient space \(H\).
This ball union is not yet useful for proving contractibility: although its
centers lie in \(X\), the balls extend away from \(X\), and their union need
not be convex.
Instead, we use a tube around \(X\) as a convex intermediate
object. We replace each ball by its intersection with this tube and show that
the restriction preserves the nerve of the cover.

\Cref{lem:ball-nerve-restriction} gives a metric condition under which
restricting the open-ball cover preserves its nerve.
\Cref{prop:finite-jung-tube-restriction} then uses the finite Jung constant
to verify this condition and obtain the required homotopy equivalence.

\begin{lemma}[Nerve-preserving open-tube restriction]
\label{lem:ball-nerve-restriction}
Let \(H\) be a hyperconvex metric space, and let \(Y\subset X\subset H\).
Fix \(\rho,\beta>0\), and suppose that, for every finite nonempty
\(W\subset Y\) with \(\operatorname{diam}(W)<2\rho\), there is \(c\in X\)
such that
\[
\max_{w\in W}d_H(c,w)<\rho+\beta.
\]
Then the covers
\[
\{\Blt{y}{H}{\rho}\}_{y\in Y}
\quad\text{and}\quad
\{\Blt{y}{H}{\rho}\cap\Tlt{X}{H}{\beta}\}_{y\in Y}
\]
have the same nerve.
\end{lemma}

\begin{proof}
Restriction to the tube $\Tlt{X}{H}{\beta}$ gives one inclusion of nerves.
For the reverse inclusion, let \(W\subset Y\) be finite and nonempty, and
suppose that
\(
\bigcap_{w\in W}\Blt{w}{H}{\rho}\neq\emptyset.
\)
Then \(\operatorname{diam}(W)<2\rho\).
Choose \(c\in X\) as in the hypothesis and
set
\[
\mu:=\max_{w\in W}d_H(c,w).
\]
Since \(\mu<\rho+\beta\), we can choose
\(
0<\varepsilon<
\min\left\{
\rho,\beta,
\rho-\frac{\operatorname{diam}(W)}2,
\frac{\rho+\beta-\mu}{2}
\right\},
\)
so that for all \(w,w'\in W\),
\[
d_H(w,w')
\le\operatorname{diam}(W)
<2(\rho-\varepsilon)
\quad\text{and}\quad
d_H(c,w)
\le\mu
<(\rho-\varepsilon)+(\beta-\varepsilon).
\]

Thus, the closed balls
\[
\{\Ble{w}{H}{\rho-\varepsilon}\}_{w\in W}
\quad\text{and}\quad
\Ble{c}{H}{\beta-\varepsilon}
\]
satisfy the pairwise hyperconvexity inequalities. Hence hyperconvexity gives a
common point \(z\in H\). Since \(\varepsilon>0\), this point $z$ lies in every
open ball \(\Blt{w}{H}{\rho}\). Moreover, since \(c\in X\) and
\(
\operatorname{dist}_H(z,X)
\le d_H(z,c)
\le\beta-\varepsilon
<\beta,
\)
$z$ also lies in
\(\Tlt{X}{H}{\beta}\). Thus the restricted balls indexed by \(W\) have
nonempty common intersection.
\end{proof}

The comparison in \cref{lem:ball-nerve-restriction} is purely metric; in
particular, \(X\) need not be a normed space. The normed structure is used in
\cref{prop:finite-jung-tube-restriction} to show that both covers are good and
hence that their unions are homotopy equivalent.

\begin{proposition}[Nerve-preserving open-tube restriction under a finite Jung bound]
\label{prop:finite-jung-tube-restriction}
Let \(H\) be a hyperconvex normed space, let \(X\subset H\) be convex, and let
\(Y\subset X\). Fix \(\rho,\beta>0\).
Suppose that every finite nonempty \(W\subset Y\) with
\(\operatorname{diam}(W)<2\rho\) satisfies
\[
R_X(W)<\rho+\beta.
\]
The circumradius hypothesis is satisfied, in particular, if
\[
J_{\mathrm{fin}}(X)<\frac{\rho+\beta}{2\rho}.
\]
Then the covers
\[
\{\Blt{y}{H}{\rho}\}_{y\in Y}
\quad\text{and}\quad
\{\Blt{y}{H}{\rho}\cap\Tlt{X}{H}{\beta}\}_{y\in Y}
\]
have the same nerve.
Moreover,
\[
\bigcup_{y\in Y}\Blt{y}{H}{\rho}
\simeq
\left(\bigcup_{y\in Y}\Blt{y}{H}{\rho}\right)
\cap\Tlt{X}{H}{\beta}.
\]
\end{proposition}

\begin{proof}
We first verify the sufficient condition involving \(J_{\mathrm{fin}}(X)\).
If
\(
J_{\mathrm{fin}}(X)<\frac{\rho+\beta}{2\rho},
\)
then every finite \(W\subset Y\) of positive diameter less than \(2\rho\)
satisfies
\[
R_X(W)
\le J_{\mathrm{fin}}(X)\operatorname{diam}(W)
<\rho+\beta.
\]
For a singleton, the same inequality follows from \(R_X(W)=0\). Thus the
bound on \(J_{\mathrm{fin}}(X)\) implies the stated circumradius hypothesis.

Assume now that the circumradius hypothesis holds. Let \(W\subset Y\) be
finite and nonempty with
\(\operatorname{diam}(W)<2\rho\). Since \(X\) carries the metric induced
from \(H\), the inequality \(R_X(W)<\rho+\beta\) gives \(c\in X\) such that
\[
\max_{w\in W}d_H(c,w)
<\rho+\beta.
\]
Thus \cref{lem:ball-nerve-restriction} shows that the two covers have the
same nerve.

The sets \(X\) and \(\Blt{0}{H}{\beta}\) are convex. Hence
\[
\Tlt{X}{H}{\beta}=X+\Blt{0}{H}{\beta}
\]
is convex as the Minkowski sum of two convex sets. It is also open
because it is a union of translates of the open ball \(\Blt{0}{H}{\beta}\).
Since the balls and the tube are convex, every nonempty finite
intersection in either cover is convex, and hence contractible.
Moreover, the balls and the tube are open, so both covers are good open covers of their unions.
These unions are metric spaces and therefore paracompact. The open
nerve theorem \cite[Corollary~4G.3]{hatcher} identifies each union, up
to homotopy, with its nerve. Since the two covers have the same nerve,
their unions are homotopy equivalent.
\end{proof}

\begin{remark}[Why a thickened tube is needed]
Directly restricting the ambient balls to the embedded subspace \(X\) can
destroy finite intersections and therefore change the nerve. The thickened
tube retains enough of the ambient space to preserve these intersections
under the hypotheses of \cref{prop:finite-jung-tube-restriction}.
The example in \cref{fig:tube-restriction-example} illustrates this
phenomenon for three balls in \((\R^3,\ell_\infty)\). Panel~\textup{(a)}
shows their common ambient intersection, panel~\textup{(b)} shows that
restriction to \(X\) destroys this intersection, and panel~\textup{(c)} shows
that restriction to a sufficiently thick tube about \(X\) preserves it. The
nerve beneath each panel records the corresponding change.
\end{remark}

\begin{figure}[htbp]
\centering
\begingroup
\definecolor{planeBlue}{RGB}{46,97,236}
\definecolor{ballRed}{RGB}{213,70,65}
\definecolor{ballOrange}{RGB}{226,126,10}
\definecolor{ballPurple}{RGB}{154,63,145}
\definecolor{tubeGreen}{RGB}{32,132,69}
\def\figureMathSize{10}
\def\panelColumnWidth{.30\textwidth}
\def\panelNerveWidth{.36\linewidth}
\def\panelNerveGap{1pt}
\def\panelBoundingBottom{3.85}
\def\panelBoundingTop{10.45}
\def\panelTitleY{10.20}
\def\panelSubtitleY{9.66}
\def\panelCardBottom{4.02}
\def\panelCardTop{9.20}
\def\panelHLabelY{9.02}
\def\panelResultY{4.25}
\def\bcOriginY{6.5cm}
\newcommand{\figuremathfont}{%
  \pgfmathsetmacro{\currentFigureMathLeading}{1.18*\figureMathSize}%
  \fontsize{\figureMathSize}{\currentFigureMathLeading}\selectfont}
\tikzset{
  panelpicture/.style={x=1cm,y=1cm,line cap=round,line join=round,>=Latex},
  paneltitle/.style={font=\sffamily\bfseries\fontsize{12.2}{14.2}\selectfont,
    anchor=west,text=black!92},
  subtitle/.style={font=\sffamily\fontsize{8.4}{10}\selectfont,
    anchor=west,text=black!57},
  objectlabel/.style={font=\figuremathfont},
  boxedobjectlabel/.style={objectlabel,fill=white,fill opacity=.78,
    text opacity=1,rounded corners=.7pt,inner sep=.8pt},
  panelcard/.style={fill=black!1,draw=black!18,rounded corners=2pt,
    line width=.65pt},
  ambientlabel/.style={font=\sffamily\figuremathfont,text=black!62,
    anchor=north west},
  panelresult/.style={font=\figuremathfont},
  coordinateaxis/.style={black!40,densely dashed,line width=.58pt,
    -{Latex[length=1.8mm]}},
}

\def\centerOne{0,0,0}
\def\centerTwo{1,0,1}
\def\centerThree{0,1,1}
\newcommand{\drawballcenters}{%
  \filldraw[fill=ballRed,draw=white,line width=.45pt]
    (\centerOne) circle[radius=.075cm];
  \filldraw[fill=ballOrange,draw=white,line width=.45pt]
    (\centerTwo) circle[radius=.075cm];
  \filldraw[fill=ballPurple,draw=white,line width=.45pt]
    (\centerThree) circle[radius=.075cm];%
}
\newcommand{\drawstandardballs}{%
  \drawcoloredcube{\centerOne}{ballRed}
  \drawcoloredcube{\centerTwo}{ballOrange}
  \drawcoloredcube{\centerThree}{ballPurple}%
}
\newcommand{\drawpanelcard}[3]{%
  \draw[panelcard]
    (#1,\panelCardBottom) rectangle (#2,\panelCardTop);
  \node[ambientlabel] at (#3,\panelHLabelY) {\(H\)};%
}
\def\bcBallOpacity{.40}
\def\aYOnePos{0.2,-.03,0}
\def\aYTwoPos{0.75,0,1.02}
\def\aYThreePos{0.21,1.03,0.95}
\def\aBOnePos{4.3,5.5}
\def\aBTwoPos{3.8,8.25}
\def\aBThreePos{2,8.5}
\def\drawcoloredcube#1#2{%
  \begin{scope}[shift={(#1)}]
    \fill[#2!72!black,fill opacity=.13]
      (-.6,-.6,-.6)--(.6,-.6,-.6)--(.6,.6,-.6)--(-.6,.6,-.6)--cycle;
    \fill[#2,fill opacity=.24]
      (-.6,-.6,.6)--(.6,-.6,.6)--(.6,.6,.6)--(-.6,.6,.6)--cycle;
    \fill[#2!88!black,fill opacity=.20]
      (-.6,-.6,-.6)--(.6,-.6,-.6)--(.6,-.6,.6)--(-.6,-.6,.6)--cycle;
    \fill[#2!78!black,fill opacity=.27]
      (.6,-.6,-.6)--(.6,.6,-.6)--(.6,.6,.6)--(.6,-.6,.6)--cycle;
    \fill[#2!55!white,fill opacity=.22]
      (.6,.6,-.6)--(-.6,.6,-.6)--(-.6,.6,.6)--(.6,.6,.6)--cycle;
    \fill[#2!68!black,fill opacity=.16]
      (-.6,.6,-.6)--(-.6,-.6,-.6)--(-.6,-.6,.6)--(-.6,.6,.6)--cycle;
    \draw[#2!88!black,line width=.62pt,opacity=.82]
      (-.6,-.6,-.6)--(.6,-.6,-.6)--(.6,.6,-.6)--(-.6,.6,-.6)--cycle
      (-.6,-.6,.6)--(.6,-.6,.6)--(.6,.6,.6)--(-.6,.6,.6)--cycle
      (-.6,-.6,-.6)--(-.6,-.6,.6)
      (.6,-.6,-.6)--(.6,-.6,.6)
      (.6,.6,-.6)--(.6,.6,.6)
      (-.6,.6,-.6)--(-.6,.6,.6);
  \end{scope}%
}
\def\drawcommoncube{%
  \begin{scope}[shift={(.5,.5,.5)}]
    \filldraw[fill=black!96,draw=black,line width=.72pt]
      (-.1,-.1,-.1)--(.1,-.1,-.1)--(.1,.1,-.1)--(-.1,.1,-.1)--cycle;
    \filldraw[fill=black!92,draw=black,line width=.72pt]
      (-.1,-.1,.1)--(.1,-.1,.1)--(.1,.1,.1)--(-.1,.1,.1)--cycle;
    \filldraw[fill=black!88,draw=black,line width=.72pt]
      (.1,-.1,-.1)--(.1,.1,-.1)--(.1,.1,.1)--(.1,-.1,.1)--cycle;
    \filldraw[fill=black!84,draw=black,line width=.72pt]
      (.1,.1,-.1)--(-.1,.1,-.1)--(-.1,.1,.1)--(.1,.1,.1)--cycle;
    \node[font=\bfseries\figuremathfont,text=white]
      at (0,0,0) {\(Q\)};
  \end{scope}%
}
\newcommand{\drawinducednerve}[1]{%
  \begin{tikzpicture}[x=1cm,y=1cm,line cap=round,line join=round]
    \path[use as bounding box] (-1.02,-.66) rectangle (1.02,.60);
    \path[#1] (-.78,.42)--(.78,.42)--(0,-.50)--cycle;
    \draw[black!88,line width=.92pt]
      (-.78,.42)--(.78,.42)--(0,-.50)--cycle;
    \filldraw[fill=ballPurple,draw=black!80,line width=.42pt]
      (-.78,.42) circle (.080);
    \filldraw[fill=ballOrange,draw=black!80,line width=.42pt]
      (.78,.42) circle (.080);
    \filldraw[fill=ballRed,draw=black!80,line width=.42pt]
      (0,-.50) circle (.080);
  \end{tikzpicture}%
}
\newcommand{\inducednervelabel}[2]{%
  \par\vspace{-2pt}%
  {\figuremathfont\textcolor{#1}{\(#2\)}\par}%
}
\newcommand{\drawpanelnerve}[3]{%
  \par\vspace{\panelNerveGap}%
  \resizebox{\panelNerveWidth}{!}{\drawinducednerve{#1}}%
  \inducednervelabel{#2}{#3}%
}
\def\bAzimuth{35}      %
\def\bElevation{52.5}     %
\def\bViewScale{1.25}
\def\bOriginX{9.25cm}
\def\bOriginY{\bcOriginY}
\def\bBOnePos{-0.32,-1.68,-.2}
\def\bBTwoPos{-0.32,-.78,.95}
\def\bBThreePos{-.78,0.98,1.6}
\def\bcXLabelPos{-1.18,1.18,-0.18}
\pgfmathsetlengthmacro{\bxx}{-\bViewScale*sin(\bAzimuth)*1cm}
\pgfmathsetlengthmacro{\bxy}{-\bViewScale*sin(\bElevation)*cos(\bAzimuth)*1cm}
\pgfmathsetlengthmacro{\byx}{ \bViewScale*cos(\bAzimuth)*1cm}
\pgfmathsetlengthmacro{\byy}{-\bViewScale*sin(\bElevation)*sin(\bAzimuth)*1cm}
\pgfmathsetlengthmacro{\bzx}{0cm}
\pgfmathsetlengthmacro{\bzy}{ \bViewScale*cos(\bElevation)*1cm}
\def\bPlaneMin{-.62}
\def\bPlaneMax{1.68}
\def\bPlaneOpacity{.15}
\pgfmathsetmacro{\bPlaneMinZ}{2*\bPlaneMin}
\pgfmathsetmacro{\bPlaneMixedZ}{\bPlaneMin+\bPlaneMax}
\pgfmathsetmacro{\bPlaneMaxZ}{2*\bPlaneMax}
\def\cViewScale{1.20}
\def\cOriginX{16.2cm}
\def\cOriginY{\bcOriginY}
\def\tubeBeta{.125}
\def\cPlaneMin{\bPlaneMin}
\def\cPlaneMax{\bPlaneMax}
\pgfmathsetmacro{\cPlaneMinZ}{2*\cPlaneMin}
\pgfmathsetmacro{\cPlaneMixedZ}{\cPlaneMin+\cPlaneMax}
\pgfmathsetmacro{\cPlaneMaxZ}{2*\cPlaneMax}
\pgfmathsetmacro{\tubeOffset}{3*\tubeBeta}
\pgfmathsetlengthmacro{\cxx}{-\cViewScale*sin(\bAzimuth)*1cm}
\pgfmathsetlengthmacro{\cxy}{-\cViewScale*sin(\bElevation)*cos(\bAzimuth)*1cm}
\pgfmathsetlengthmacro{\cyx}{ \cViewScale*cos(\bAzimuth)*1cm}
\pgfmathsetlengthmacro{\cyy}{-\cViewScale*sin(\bElevation)*sin(\bAzimuth)*1cm}
\pgfmathsetlengthmacro{\czx}{0cm}
\pgfmathsetlengthmacro{\czy}{ \cViewScale*cos(\bElevation)*1cm}
\pgfmathsetmacro{\cLowerMinZ}{\cPlaneMinZ-\tubeOffset}
\pgfmathsetmacro{\cLowerMixedZ}{\cPlaneMixedZ-\tubeOffset}
\pgfmathsetmacro{\cLowerMaxZ}{\cPlaneMaxZ-\tubeOffset}
\pgfmathsetmacro{\cUpperMinZ}{\cPlaneMinZ+\tubeOffset}
\pgfmathsetmacro{\cUpperMixedZ}{\cPlaneMixedZ+\tubeOffset}
\pgfmathsetmacro{\cUpperMaxZ}{\cPlaneMaxZ+\tubeOffset}

\def\aPanelLeft{.12}
\def\aPanelRight{5.72}
\def\bPanelLeft{6.67}
\def\bPanelRight{12.27}
\def\cPanelLeft{13.42}
\def\cPanelRight{19.02}

\noindent
\begin{minipage}[t]{\panelColumnWidth}
  \vspace{0pt}\centering
  \resizebox{\linewidth}{!}{%
    \begin{tikzpicture}[panelpicture]
    \path[use as bounding box]
      (0,\panelBoundingBottom) rectangle (5.84,\panelBoundingTop);
    \node[paneltitle] at (.16,\panelTitleY) {(a)\; Ambient ball cover};
    \node[subtitle] at (.19,\panelSubtitleY)
      {the triple intersection is nonempty};
    \drawpanelcard{\aPanelLeft}{\aPanelRight}{.32}

    \begin{scope}[shift={(2.58,5.80)},
        x={(1.10cm,.15cm)},y={(-.65cm,.40cm)},z={(0cm,1.10cm)}]
      \draw[coordinateaxis]
        (-1.05,0,0)--(2.05,0,0);
      \draw[coordinateaxis]
        (0,-1.05,0)--(0,2.05,0);
      \draw[coordinateaxis]
        (0,0,-.90)--(0,0,2.05);

      \drawcoloredcube{\centerOne}{ballRed}
      \drawcoloredcube{\centerThree}{ballPurple}
      \drawcoloredcube{\centerTwo}{ballOrange}

      \drawcommoncube

      \drawballcenters
      \node[objectlabel,text=ballRed,anchor=north east]
        at (\aYOnePos) {\(y_1\)};
      \node[objectlabel,text=ballOrange,anchor=south west]
        at (\aYTwoPos) {\(y_2\)};
      \node[objectlabel,text=ballPurple,anchor=south east]
        at (\aYThreePos) {\(y_3\)};
    \end{scope}
    \node[objectlabel,text=ballRed!92!black,anchor=north east]
      at (\aBOnePos) {\(B_1\)};
    \node[objectlabel,text=ballOrange!92!black,anchor=west]
      at (\aBTwoPos) {\(B_2\)};
    \node[objectlabel,text=ballPurple!92!black,anchor=east]
      at (\aBThreePos) {\(B_3\)};
    \node[panelresult,text=black!74] at (2.92,\panelResultY)
      {\(\varnothing\neq Q=
        \textcolor{ballRed}{B_1}\cap
        \textcolor{ballOrange}{B_2}\cap
        \textcolor{ballPurple}{B_3}\)};
    \end{tikzpicture}
  }
  \drawpanelnerve{fill=black!18}{black!75}{\mathcal N(\{B_i\})}
\end{minipage}%
\hfill
\begin{minipage}[t]{\panelColumnWidth}
  \vspace{0pt}\centering
  \resizebox{\linewidth}{!}{%
    \begin{tikzpicture}[panelpicture]
    \path[use as bounding box]
      (6.55,\panelBoundingBottom) rectangle (12.39,\panelBoundingTop);
    \node[paneltitle] at (6.71,\panelTitleY) {(b)\; Restriction to \(X\)};
    \node[subtitle] at (6.74,\panelSubtitleY)
      {the triple intersection is lost};
    \drawpanelcard{\bPanelLeft}{\bPanelRight}{6.87}

    \begin{scope}[shift={(\bOriginX,\bOriginY)},
        x={(\bxx,\bxy)},y={(\byx,\byy)},z={(\bzx,\bzy)}]
      \draw[coordinateaxis]
        (-1.00,0,0)--(2.00,0,0);
      \draw[coordinateaxis]
        (0,-1.00,0)--(0,2.00,0);
      \draw[coordinateaxis]
        (0,0,-.78)--(0,0,2.05);

      \begin{scope}[transparency group,opacity=\bcBallOpacity]
        \drawstandardballs
      \end{scope}

      \filldraw[fill=planeBlue,fill opacity=\bPlaneOpacity,draw=planeBlue!82,
        line width=.84pt]
        (\bPlaneMin,\bPlaneMin,\bPlaneMinZ)--
        (\bPlaneMax,\bPlaneMin,\bPlaneMixedZ)--
        (\bPlaneMax,\bPlaneMax,\bPlaneMaxZ)--
        (\bPlaneMin,\bPlaneMax,\bPlaneMixedZ)--cycle;

      \drawcommoncube

      \drawballcenters
      \node[boxedobjectlabel,text=ballRed!92!black,anchor=north east]
        at (\bBOnePos) {\(B_1\)};
      \node[boxedobjectlabel,text=ballOrange!92!black,anchor=south west]
        at (\bBTwoPos) {\(B_2\)};
      \node[boxedobjectlabel,text=ballPurple!92!black,anchor=south west]
        at (\bBThreePos) {\(B_3\)};
      \node[objectlabel,text=planeBlue,anchor=south west]
        at (\bcXLabelPos) {\(X\)};
    \end{scope}
    \node[panelresult,text=black!74] at (9.47,\panelResultY)
      {\(Q\cap\textcolor{planeBlue}{X}=\varnothing\)};
    \end{tikzpicture}
  }
  \drawpanelnerve{fill=none}{planeBlue}{\mathcal N(\{B_i\cap X\})}
\end{minipage}%
\hfill
\begin{minipage}[t]{\panelColumnWidth}
  \vspace{0pt}\centering
  \resizebox{\linewidth}{!}{%
    \begin{tikzpicture}[panelpicture]
    \path[use as bounding box]
      (13.30,\panelBoundingBottom) rectangle (19.14,\panelBoundingTop);
    \node[paneltitle] at (13.46,\panelTitleY)
      {(c)\; Restriction to \(T_\beta(X)\)};
    \node[subtitle] at (13.49,\panelSubtitleY)
      {the triple intersection is preserved};
    \drawpanelcard{\cPanelLeft}{\cPanelRight}{13.64}

    \begin{scope}[shift={(\cOriginX,\cOriginY)},
        x={(\cxx,\cxy)},y={(\cyx,\cyy)},z={(\czx,\czy)}]
      \drawcommoncube

      \fill[tubeGreen,fill opacity=.12]
        (\cPlaneMin,\cPlaneMin,\cLowerMinZ)--
        (\cPlaneMax,\cPlaneMin,\cLowerMixedZ)--
        (\cPlaneMax,\cPlaneMin,\cUpperMixedZ)--
        (\cPlaneMin,\cPlaneMin,\cUpperMinZ)--cycle;
      \fill[tubeGreen,fill opacity=.15]
        (\cPlaneMax,\cPlaneMin,\cLowerMixedZ)--
        (\cPlaneMax,\cPlaneMax,\cLowerMaxZ)--
        (\cPlaneMax,\cPlaneMax,\cUpperMaxZ)--
        (\cPlaneMax,\cPlaneMin,\cUpperMixedZ)--cycle;
      \fill[tubeGreen,fill opacity=.12]
        (\cPlaneMin,\cPlaneMax,\cLowerMixedZ)--
        (\cPlaneMax,\cPlaneMax,\cLowerMaxZ)--
        (\cPlaneMax,\cPlaneMax,\cUpperMaxZ)--
        (\cPlaneMin,\cPlaneMax,\cUpperMixedZ)--cycle;
      \fill[tubeGreen,fill opacity=.15]
        (\cPlaneMin,\cPlaneMin,\cLowerMinZ)--
        (\cPlaneMin,\cPlaneMax,\cLowerMixedZ)--
        (\cPlaneMin,\cPlaneMax,\cUpperMixedZ)--
        (\cPlaneMin,\cPlaneMin,\cUpperMinZ)--cycle;
      \filldraw[fill=tubeGreen,fill opacity=.12,draw=tubeGreen!78!black,
        line width=.72pt]
        (\cPlaneMin,\cPlaneMin,\cLowerMinZ)--
        (\cPlaneMax,\cPlaneMin,\cLowerMixedZ)--
        (\cPlaneMax,\cPlaneMax,\cLowerMaxZ)--
        (\cPlaneMin,\cPlaneMax,\cLowerMixedZ)--cycle;
      \filldraw[fill=tubeGreen,fill opacity=.12,draw=tubeGreen!78!black,
        line width=.72pt]
        (\cPlaneMin,\cPlaneMin,\cUpperMinZ)--
        (\cPlaneMax,\cPlaneMin,\cUpperMixedZ)--
        (\cPlaneMax,\cPlaneMax,\cUpperMaxZ)--
        (\cPlaneMin,\cPlaneMax,\cUpperMixedZ)--cycle;
      \draw[tubeGreen!72!black,densely dotted,line width=.68pt]
        (\cPlaneMin,\cPlaneMin,\cLowerMinZ)--
        (\cPlaneMin,\cPlaneMin,\cUpperMinZ);
      \draw[tubeGreen!72!black,densely dotted,line width=.68pt]
        (\cPlaneMax,\cPlaneMin,\cLowerMixedZ)--
        (\cPlaneMax,\cPlaneMin,\cUpperMixedZ);
      \draw[tubeGreen!72!black,densely dotted,line width=.68pt]
        (\cPlaneMax,\cPlaneMax,\cLowerMaxZ)--
        (\cPlaneMax,\cPlaneMax,\cUpperMaxZ);
      \draw[tubeGreen!72!black,densely dotted,line width=.68pt]
        (\cPlaneMin,\cPlaneMax,\cLowerMixedZ)--
        (\cPlaneMin,\cPlaneMax,\cUpperMixedZ);

      \begin{scope}[transparency group,opacity=\bcBallOpacity]
        \drawstandardballs
      \end{scope}

      \filldraw[fill=planeBlue,fill opacity=.16,draw=planeBlue!88,
        densely dashed,line width=1.05pt]
        (\cPlaneMin,\cPlaneMin,\cPlaneMinZ)--
        (\cPlaneMax,\cPlaneMin,\cPlaneMixedZ)--
        (\cPlaneMax,\cPlaneMax,\cPlaneMaxZ)--
        (\cPlaneMin,\cPlaneMax,\cPlaneMixedZ)--cycle;
      \drawballcenters
      \node[objectlabel,text=planeBlue,anchor=south west]
        at (\bcXLabelPos) {\(X\)};
      \node[objectlabel,text=tubeGreen!82!black,anchor=north east,
        fill=white,fill opacity=.78,text opacity=1,inner sep=.8pt]
        at (5.42,2.08,8.05) {\(T_\beta(X)\)};
    \end{scope}
    \node[panelresult] at (16.22,\panelResultY)
      {\(
        \textcolor{tubeGreen!82!black}{T_{\beta}(X)}
        \cap Q \neq \varnothing \)};
    \end{tikzpicture}
  }
  \drawpanelnerve{fill=tubeGreen!24}{tubeGreen!84!black}
    {\mathcal N(\{B_i\cap T_\beta(X)\})}
\end{minipage}
\endgroup

\caption{{\textbf{Restriction to a tube can preserve the nerve when restriction to the
embedded subspace does not.}} Let \(H=(\R^3,\ell_\infty)\), let
\(X=\{(u,v,u+v):u,v\in\R\}\), and set \(y_1=(0,0,0)\), \(y_2=(1,0,1)\),
and \(y_3=(0,1,1)\), with \(B_i:=\Blt{y_i}{H}{3/5}\). In each panel, the
upper row shows the relevant cover and the lower row shows its nerve.
\textup{(a)} The ambient balls have common intersection
\(Q=(2/5,3/5)^3\), represented by the dark cube, so their nerve contains a
\(2\)-simplex. \textup{(b)} Since \(Q\cap X=\varnothing\), restriction to
\(X\) destroys the triple intersection and removes the \(2\)-simplex.
\textup{(c)} Set \(\beta=1/8\). The point \(q=(5/12,5/12,7/12)\) belongs to
\(Q\) and satisfies \(\operatorname{dist}_{\ell_\infty}(q,X)=1/12<\beta\).
Hence the triple intersection survives after restriction to \(T_\beta(X)\),
and the \(2\)-simplex is preserved. Boundaries are drawn only for visibility.}
\label{fig:tube-restriction-example}
\end{figure}

\subsection{Density fills the open tube}
\label{subsec:open-tube-density}

We now combine the tube restriction with the \(\delta\)-density of \(Y\) to prove \cref{thm:jung-open-rips}. The tube must be chosen large enough to preserve the nerve, but small enough that the \(\rho\)-balls centered at \(Y\) cover it. We then apply Bohnenblust's finite-dimensional Jung bound to obtain an explicit scale depending only on the dimension and density.

\openjungcontractibility*

\begin{proof}%
If \(X=\{0\}\), then \(\delta\)-density forces \(Y=\{0\}\), and the result is clear.
Assume \(X\neq\{0\}\). Then
\(J_{\mathrm{fin}}(X)\ge1/2\), because every two-point set has circumradius
equal to half its diameter. The scale hypothesis allows us to choose \(J\)
such that
\[
J_{\mathrm{fin}}(X)<J<1-\frac{\delta}{r}.
\]
Set
\[
\rho:=\frac r2,
\qquad
\beta:=(2J-1)\rho .
\]
Then \(\beta>0\), \(\delta+\beta<\rho\), and
\[
J=\frac{\rho+\beta}{2\rho}.
\]

Let \(H:=\ell_\infty(B_{X^*})\), and let
\(\iota:X\to H\) be the canonical isometric embedding (see \cref{lem:canonical-embedding}). Write
\[
T:=\Tlt{\iota(X)}{H}{\beta}.
\]
The image \(\iota(X)\) is a linear subspace of \(H\), hence convex.
Since \(\Blt{0}{H}{\beta}\) is also convex, their Minkowski sum
\(
T=\iota(X)+\Blt{0}{H}{\beta}
\)
is convex.

The space \(H\) is hyperconvex by
\cite[Example~2.22 and Proposition~2.18]{limmemoliokutan}.
Let \(W\subset\iota(Y)\) be finite and nonempty with
\(\operatorname{diam}(W)<2\rho\). If \(W\) has positive diameter, then
\[
R_{\iota(X)}(W)
\le J_{\mathrm{fin}}(\iota(X))\operatorname{diam}(W)
<2J\rho
=\rho+\beta.
\]
For a singleton, the same inequality follows from
\(R_{\iota(X)}(W)=0\). Thus
\cref{prop:finite-jung-tube-restriction}, applied to
\(\iota(Y)\subset\iota(X)\subset H\), and
\cref{prop:open-thickening} give
\[
\VRlt{Y}{2\rho}
\simeq \bigcup_{y\in Y}\Blt{\iota(y)}{H}{\rho}
\simeq
\left(\bigcup_{y\in Y}\Blt{\iota(y)}{H}{\rho}\right)\cap T.
\]

It remains to identify the restricted union with \(T\). Since the restricted
union is contained in \(T\) by definition, it suffices to prove
\[
T\subset\bigcup_{y\in Y}\Blt{\iota(y)}{H}{\rho}.
\]
Let \(z\in T\).
Choose \(x\in X\) with \(d_H(z,\iota(x))<\beta\), and then choose
\(y\in Y\) with \(d_X(x,y)\le\delta\). Since \(\iota\) is an isometry and
\(\delta+\beta<\rho\),
\[
\begin{aligned}
    d_H(z,\iota(y))
     & \le d_H(z,\iota(x))+d_H(\iota(x),\iota(y)) \\
     & =d_H(z,\iota(x))+d_X(x,y)                  \\
     & <\beta+\delta<\rho.
\end{aligned}
\]
Thus \(z\in\Blt{\iota(y)}{H}{\rho}\), and hence
\[
\left(\bigcup_{y\in Y}\Blt{\iota(y)}{H}{\rho}\right)\cap T=T.
\]
The space \(T\) is convex and therefore contractible. Since \(r=2\rho\), the
result follows.
\end{proof}

\begin{remark}
The proof of \cref{thm:jung-open-rips} applies without change to any normed
space \(X\) satisfying \(J_{\mathrm{fin}}(X)<1\). We state the theorem in
finite-dimensional form because the applications considered in this paper are
finite-dimensional. In this setting, Bohnenblust's bound also provides the
explicit dimension-dependent estimate used below.
\end{remark}

To obtain a bound depending only on dimension, we now invoke the classical
Jung inequality for normed spaces.

\begin{lemma}[Bohnenblust's finite-dimensional Jung bound
{\cite[Theorem~6]{bohnenblust}}]
\label{lem:bohnenblust-general}
Let \(X\) be an \(n\)-dimensional normed space with \(n\ge1\). For every
finite nonempty \(W\subset X\), there is \(c\in X\) such that
\[
\max_{w\in W}\|c-w\|
\le\frac{n}{n+1}\operatorname{diam}(W).
\]
In particular,
\[
J_{\mathrm{fin}}(X)\le\frac{n}{n+1}.
\]
\end{lemma}

\begin{corollary}[{Finite-dimensional bound for open
Vietoris--Rips complexes}]
\label{cor:finite-dimensional-open}
Let \(X\) be an \(n\)-dimensional normed space with \(n\ge1\), and let
\(Y\subset X\) be \(\delta\)-dense. Then \(\VRlt{Y}{r}\) is contractible for every
\[
r>(n+1)\delta.
\]
\end{corollary}

\begin{proof}
By \cref{lem:bohnenblust-general},
\[
J_{\mathrm{fin}}(X)\le\frac{n}{n+1}.
\]
Thus \(\delta/(1-J_{\mathrm{fin}}(X))\le(n+1)\delta\), and the result follows from
\cref{thm:jung-open-rips}.
\end{proof}

\subsection{\texorpdfstring{Explicit estimates in
            \(\ell_p^n\)}{Explicit estimates in lp spaces}}
\label{subsec:lp-open-estimates}

Explicit applications of \cref{thm:jung-open-rips} to \(\delta\)-dense subsets
of \(\ell_p^n\) require estimates for \(J_{\mathrm{fin}}(\ell_p^n)\). The Jung
inequality of Ivanov and Pichugov gives one such estimate, while the
inequalities comparing the \(\ell_p\)- and \(\ell_\infty\)-norms give another.
Their minimum yields the general bound in \cref{cor:lp-open-dense}. In the
planar \(\ell_1\) case, an isometry with \(\ell_\infty^2\) gives the exact
finite Jung constant and improves this general estimate.

\begin{lemma}[Jung estimates for \(\ell_p^n\)]
\label{lem:lp-jung}

Let \(n\ge1\) and \(1\le p<\infty\). Define \(p'\in(1,\infty]\) by
\(1/p+1/p'=1\), where \(1/\infty:=0\), and set
\(s_p:=\max\{p,p'\}\). Define
\[
j_{p,n}:=\min\left\{
2^{-1/s_p}\left(\frac{n}{n+1}\right)^{1/p},
\frac{n^{1/p}}2
\right\}.
\]
Then
\[
J_{\mathrm{fin}}(\ell_p^n)\le j_{p,n}.
\]
For \(p=\infty\), one has
\(J_{\mathrm{fin}}(\ell_\infty^n)=1/2\).

\end{lemma}

\begin{proof}
It suffices to bound \(J_{\mathrm{fin}}(\ell_p^n)\) by each term in the
definition of \(j_{p,n}\). The first bound follows from the Jung inequality
of Ivanov and Pichugov~\cite[inequalities~(3)--(4)]{ivanovpichugov}. For the
second, let \(W\subset\R^n\) be finite with positive diameter.
Hyperconvexity of \(\ell_\infty^n\) gives \(c\in\R^n\) such that
\[
\max_{w\in W}\|c-w\|_\infty
\le\frac12\operatorname{diam}_\infty(W).
\]
Since \(\|x\|_p\le n^{1/p}\|x\|_\infty\) and
\(\operatorname{diam}_\infty(W)\le\operatorname{diam}_p(W)\),
\[
\max_{w\in W}\|c-w\|_p
\le\frac{n^{1/p}}2\operatorname{diam}_p(W).
\]
Taking the supremum over all such \(W\) proves the assertion for finite
\(p\). For \(p=\infty\), the same hyperconvexity argument gives the upper
bound \(1/2\), while any two-point set gives the reverse inequality.
\end{proof}

\begin{corollary}[Open bounds in \(\ell_p^n\)]
\label{cor:lp-open-dense}
Let \(Y\subset\ell_p^n\) be \(\delta\)-dense. If
\(1\le p<\infty\), then \(\VRlt{Y,\ell_p}{r}\) is contractible whenever
\[
r>\frac{\delta}{1-j_{p,n}}.
\]
For \(p=\infty\), the complex \(\VRlt{Y,\ell_\infty}{r}\) is contractible
whenever \(r>2\delta\).
\end{corollary}

\begin{proof}
For finite \(p\), the first term in the definition of \(j_{p,n}\) is less
than \(1\), so \(j_{p,n}<1\). The conclusion follows from
\cref{thm:jung-open-rips,lem:lp-jung}. For \(p=\infty\), use
\(J_{\mathrm{fin}}(\ell_\infty^n)=1/2\) in
\cref{thm:jung-open-rips}.
\end{proof}

\begin{lemma}[The two-dimensional sign-coordinate isometry]
\label{lem:herrlich-embedding}
The map
\[
A_2:\ell_1^2\longrightarrow \ell_\infty^2,
\qquad
A_2(x_1,x_2)=(x_1+x_2,x_1-x_2),
\]
is a surjective linear isometry.
\end{lemma}

\begin{proof}
For every \(x\in\R^2\),
\[
\|A_2x\|_\infty
=\max\{|x_1+x_2|,|x_1-x_2|\}
=|x_1|+|x_2|.
\]
The matrix of \(A_2\) has nonzero determinant, so \(A_2\) is surjective.
\end{proof}

The preceding isometry improves the general estimate in
\cref{cor:lp-open-dense} for planar \(\ell_1\).

\begin{corollary}[Open bound in \(\ell_1^2\)]
\label{cor:open-l12-sharp}
Let \(Y\subset \ell_1^2\) be \(\delta\)-dense. Then
\(\VRlt{Y}{r}\) is contractible for every \(r>2\delta\).
\end{corollary}

\begin{proof}
By \cref{lem:herrlich-embedding}, \(\ell_1^2\) and \(\ell_\infty^2\) are
isometric and therefore have the same finite Jung constant. Hence
\cref{lem:lp-jung} gives
\(
J_{\mathrm{fin}}(\ell_1^2)
=
J_{\mathrm{fin}}(\ell_\infty^2)
=
\frac12.
\)
The conclusion now follows from \cref{thm:jung-open-rips}.
\end{proof}

\section{Lattices in Normed Spaces}
\label{sec:lp-lattice-consequences}
\label{sec:l1-lattice-consequences}

Every full-rank affine lattice in a finite-dimensional normed space is
\(\delta\)-dense for some \(\delta\ge0\) and has a locally finite distance
spectrum. The density condition gives a strict open bound from
\cref{thm:jung-open-rips}.
The simplex condition in a Vietoris--Rips complex
changes only when the scale crosses a realized distance, so local finiteness
of the distance spectrum converts this strict open bound into the corresponding
weak closed bound. The following proposition isolates this observation for
arbitrary metric spaces.
For coordinate lattices, coordinatewise rounding makes the density bound
explicit (\cref{subsec:coordinate-lattice-bounds}), and specific gaps between
realized distances can improve the resulting scales
(\cref{subsec:distance-gap-bounds}). In dimensions one and two, we determine
the exact thresholds; see \cref{subsec:exact-lattice-thresholds}.

\begin{proposition}[Open and closed scales for locally finite distance spectra]
\label{prop:open-closed-threshold-equivalence}
Let \((M,d_M)\) be a metric space, and set
\(
\Sigma(M):=\{d_M(x,y):x,y\in M\}.
\)
Suppose that \(\Sigma(M)\cap[0,T]\) is finite for every \(T<\infty\).
For every \(0<a<b\), there exist
\[
r\in(a,b]
\qquad\text{and}\qquad
s\in[a,b)
\]
such that
\[
\VRle{M}{a}=\VRlt{M}{r}
\qquad\text{and}\qquad
\VRlt{M}{b}=\VRle{M}{s}.
\]
\end{proposition}

\begin{proof}
If \(\Sigma(M)\cap(a,b]\) is nonempty, let \(r\) be its least element;
otherwise, set \(r=b\). Then no distance in \(M\) lies in \((a,r)\), and
therefore
\[
\VRle{M}{a}=\VRlt{M}{r}.
\]

Similarly, if \(\Sigma(M)\cap[a,b)\) is nonempty, let \(s\) be its greatest
element; otherwise, set \(s=a\). Then no distance in \(M\) lies in
\((s,b)\), and therefore
\[
\VRlt{M}{b}=\VRle{M}{s}. \qedhere
\]
\end{proof}

For a \(\delta\)-dense full-rank affine lattice \(\Lambda\subset X\),
\cref{thm:jung-open-rips} gives the strict open bound, and
\cref{prop:open-closed-threshold-equivalence} gives the same numerical bound
with a weak inequality for the closed complex.

\normedlatticecontractibility*

\begin{proof}%
{
If \(X=\{0\}\), then \(\Lambda\) is a singleton, so both
\(\VRlt{\Lambda}{r}\) and \(\VRle{\Lambda}{r}\) are simplices, and hence
contractible, for every \(r>0\). 
Otherwise, \cref{thm:jung-open-rips} gives contractibility of
\(\VRlt{\Lambda}{r}\) for
\[
r>\frac{\delta}{1-J_{\mathrm{fin}}(X)}.
\]
The difference lattice \(\Lambda-\Lambda\) contains only finitely many
vectors in every bounded set, so the distance spectrum of \(\Lambda\) satisfies
the hypothesis of \cref{prop:open-closed-threshold-equivalence}. That
proposition gives the closed assertion.
}
\end{proof}

\subsection{\texorpdfstring{Density and bounds for
            coordinate lattices}{Density and bounds for coordinate lattices}}
\label{subsec:coordinate-lattice-bounds}

For coordinate lattices, coordinatewise rounding makes the density parameter
in the preceding contractibility results explicit. Combining this estimate
with the open bounds and the distance-spectrum comparison yields corresponding
open and closed bounds for both \(\Z^n\) and rectangular coordinate lattices.

\begin{lemma}[Density of rectangular coordinate lattices]
\label{lem:coordinate-lattice-density}
Let \(n\ge1\), \(a\in\R^n\), and \(h=(h_1,\dots,h_n)\in \R^n_{>0}\). Set
\(
\Lambda:=a+\prod_{i=1}^n h_i\Z.
\)
For \(1\le p\leq\infty\), the lattice \(\Lambda\) is \(\delta_p\)-dense in
\(\ell_p^n\), where
\(
\delta_p:=\frac12 \|h\|_p.
\)
\end{lemma}

\begin{proof}
Given \(x\in\R^n\), choose \(\lambda_i\in a_i+h_i\Z\) with
\(|x_i-\lambda_i|\le h_i/2\), and set
\(\lambda=(\lambda_1,\dots,\lambda_n)\in\Lambda\). Then
\[
\|x-\lambda\|_p
\le\frac12\left(\sum_{i=1}^n h_i^p\right)^{1/p}
\]
for finite \(p\), and
\[
\|x-\lambda\|_\infty\le\frac12\max_i h_i. \qedhere
\]
\end{proof}

Taking \(a=0\) and \(h_i=1\) shows that \(\Z^n\) is
\(n^{1/p}/2\)-dense for finite \(p\) and \(1/2\)-dense for \(p=\infty\).
Together with
\cref{cor:lp-open-dense,prop:open-closed-threshold-equivalence}, these values
give the following bounds for the integer lattice.

\lplatticecontractibility*

\begin{proof}%
For finite \(p\), \cref{lem:coordinate-lattice-density} and
\cref{cor:lp-open-dense} give contractibility of the open complex for
\(r>r_{p,n}\). The distance spectrum of \(\Z^n\) is locally finite, so
\cref{prop:open-closed-threshold-equivalence} gives the closed assertion.

For \(p=\infty\), \cref{lem:coordinate-lattice-density,cor:lp-open-dense}
give contractibility of the open complex for \(r>1\), and
\cref{prop:open-closed-threshold-equivalence} gives contractibility of the
closed complex for \(r\ge1\).
Distinct points of \(\Z^n\) have \(\ell_\infty\)-distance at least \(1\).
Thus the closed complex is discrete for \(r<1\), and the open complex is
discrete for \(r\le1\), proving necessity in both conventions.
\end{proof}

For fixed \(n\), the two terms defining \(j_{p,n}\) converge to \(1\) and
\(1/2\), respectively. Hence
\[
    j_{p,n}\longrightarrow\frac12
    \qquad\text{and}\qquad
    r_{p,n}\longrightarrow1.
\]
Thus the sufficient scales \(r_{p,n}\) converge to the exact
\(\ell_\infty\) threshold \(1\).

The bounds for \(\Z^n\) extend to rectangular coordinate lattices with
unequal coordinate spacings. Applying
\cref{cor:lp-open-dense,prop:open-closed-threshold-equivalence} with the
density values from \cref{lem:coordinate-lattice-density} gives the following
result.

\begin{corollary}[Rectangular coordinate lattices]
\label{rem:rectangular-lattices}
\label{cor:rectangular-lattice}
Let \(n\ge1\), let \(a\in\R^n\), and let \(h_1,\dots,h_n>0\). Set
\[
\Lambda=a+\prod_{i=1}^n h_i\Z\subset\R^n.
\]
For \(1\le p<\infty\), the complex
\(\VRlt{\Lambda,\ell_p}{r}\) is contractible when \(r\) is strictly larger than
\[
\frac{\left(\sum_{i=1}^n h_i^p\right)^{1/p}}
     {2(1-j_{p,n})},
\]
and \(\VRle{\Lambda,\ell_p}{r}\) is contractible when \(r\) is at least this
value. Here \(j_{p,n}\) is defined in \cref{lem:lp-jung}. For \(p=\infty\),
\(\VRlt{\Lambda,\ell_\infty}{r}\) is contractible if and only if
\(r>\max_i h_i\), while \(\VRle{\Lambda,\ell_\infty}{r}\) is contractible if
and only if \(r\ge\max_i h_i\).
\end{corollary}

    \begin{proof}
For finite \(p\), \cref{lem:coordinate-lattice-density} and
\cref{cor:lp-open-dense} give the stated open bound. The distance spectrum of
\(\Lambda\) is locally finite, so
\cref{prop:open-closed-threshold-equivalence} gives the closed bound.

For \(p=\infty\), the same argument gives contractibility of the open complex
for \(r>\max_i h_i\) and of the closed complex for \(r\ge\max_i h_i\).
Conversely, choose \(j\) such that
\(h_j=\max_i h_i\). If two lattice points have different \(j\)-th
coordinates, then their \(\ell_\infty\)-distance is at least \(h_j\).
Hence the closed complex is disconnected for \(r<h_j\), and the open complex
is disconnected for \(r\le h_j\). This proves necessity in both conventions.
\end{proof}

For \(p=1\), the threshold in \cref{cor:rectangular-lattice} simplifies to
\[
\frac{n+1}{2}\sum_{i=1}^n h_i,
\]
with a strict inequality in the open convention and a weak inequality in the
closed convention.
For fixed \(h_1,\dots,h_n\), the finite-\(p\) bound in
\cref{cor:rectangular-lattice} converges to \(\max_i h_i\) as
\(p\to\infty\), agreeing with the exact open and closed thresholds for the
rectangular coordinate lattice in the \(\ell_\infty\) metric.

\subsection{\texorpdfstring{Bounds from distance gaps}
{Bounds from distance gaps}}
\label{subsec:distance-gap-bounds}

In this subsection, we improve the sufficient scales \(r_{p,n}\) in
\cref{thm:lp-lattice} when \(p\) is a positive integer. The improvement uses
gaps in the lattice distance spectrum to sharpen the circumradius estimate. We
first formulate the underlying criterion for a general dense subset of a
finite-dimensional normed space.

Let \(Y\subset X\) be \(\delta\)-dense, and suppose that no pairwise distance
in \(Y\) lies in \((a,b)\), where \(0<a<b\). Then
\[
\VRle{Y}{a}=\VRlt{Y}{b}
\qquad(a<b).
\]
This identity uses only the single gap \((a,b)\); it does not require the
full distance spectrum to be locally finite.
Since
\(\VRlt{Y}{b}=\VRle{Y}{a}\), every simplex represented by the open
\(b/2\)-ball cover has diameter at most \(a\).
The finite Jung estimate therefore
bounds the circumradius of each simplex by
\(J_{\mathrm{fin}}(X)a\), and the tube argument requires
\[
J_{\mathrm{fin}}(X)a+\delta<b.
\]
This condition can hold even when
\(J_{\mathrm{fin}}(X)b+\delta<b\) fails. If, in addition,
\(2\delta<b\), then a positive tube width can be chosen, and the tube
argument can establish contractibility of \(\VRlt{Y}{b}\) even when
\cref{thm:jung-open-rips} does not apply at scale \(b\).

\begin{proposition}[Contractibility from a distance gap]
\label{prop:distance-gap-contractibility}
Let \(X\) be a finite-dimensional normed space, let \(Y\subset X\) be
\(\delta\)-dense, and let \(0<a<b\). Suppose that no pairwise distance in
\(Y\) belongs to \((a,b)\). If
\[
2\delta<b
\qquad\text{and}\qquad
J_{\mathrm{fin}}(X)a+\delta<b,
\]
then
\[
\VRle{Y}{a}=\VRlt{Y}{b}\simeq *.
\]
Consequently, \(\VRle{Y}{r}\) is contractible for every \(a\le r<b\), and
\(\VRlt{Y}{s}\) is contractible for every \(a<s\le b\).
\end{proposition}

\begin{proof}
Put \(J:=J_{\mathrm{fin}}(X)\) and \(\rho:=b/2\).
For every finite nonempty \(W\subset Y\) with
\(\operatorname{diam}(W)<b\), the distance-gap hypothesis gives
\(\operatorname{diam}(W)\le a\). If \(W\) has positive diameter, then
\[
R_X(W)
\le J\operatorname{diam}(W)
\le Ja.
\]
For a singleton, the same bound follows from \(R_X(W)=0\). Thus
\(R_X(W)\le Ja\) in all cases.

To apply \cref{prop:finite-jung-tube-restriction}, it is enough to choose
\(\beta>0\) such that \(Ja<\rho+\beta\). We also require this $\beta$ to satisfy
\(\delta+\beta<\rho\), which ensures that the open \(\rho\)-balls centered
at \(Y\) cover the open \(\beta\)-tube around \(X\).
Since \(2\delta<b=2\rho\) and \(Ja+\delta<b=2\rho\), both \(0\) and
\(Ja-\rho\) are strictly smaller than \(\rho-\delta\). Choose
\[
\max\{0,Ja-\rho\}<\beta<\rho-\delta.
\]
Then
\[
Ja<\rho+\beta
\qquad\text{and}\qquad
\delta+\beta<\rho.
\]

Identify \(X\) with its canonical image in the hyperconvex space
\(H:=\ell_\infty(B_{X^*})\). Since the canonical embedding is linear, its
image is a linear subspace of \(H\), and hence convex. Set
\(
T:=\Tlt{X}{H}{\beta}.
\)
Since \(H\) is hyperconvex, \cref{prop:open-thickening} gives the first
homotopy equivalence below. Since \(X\) is convex and
\(R_X(W)\le Ja<\rho+\beta\) whenever
\(\operatorname{diam}(W)<b=2\rho\),
\cref{prop:finite-jung-tube-restriction} gives the second:
\[
\VRlt{Y}{b}
\simeq
\bigcup_{y\in Y}\Blt{y}{H}{\rho}
\simeq
\left(\bigcup_{y\in Y}\Blt{y}{H}{\rho}\right)\cap T.
\]

We next show that the restricted union equals \(T\). Let \(z\in T\). Choose \(x\in X\)
with \(d_H(z,x)<\beta\), and then choose \(y\in Y\) with
\(d_H(x,y)\le\delta\). Since \(\delta+\beta<\rho\),
\[
d_H(z,y)<\beta+\delta<\rho.
\]
Thus \(z\in\Blt{y}{H}{\rho}\), so the restricted union is \(T\).
Since \(X\) and \(\Blt{0}{H}{\beta}\) are convex,
\(
T=X+\Blt{0}{H}{\beta}
\)
is convex and hence contractible.
Therefore
\(\VRlt{Y}{b}\simeq *\).

The distance gap gives
\(\VRle{Y}{a}=\VRlt{Y}{b}\). It also gives
\[
\VRle{Y}{r}=\VRle{Y}{a}\quad(a\le r<b),
\qquad
\VRlt{Y}{s}=\VRlt{Y}{b}\quad(a<s\le b).
\]
The remaining assertions follow.
\end{proof}

Since \cref{prop:distance-gap-contractibility} treats one distance gap at a
time, proving contractibility for \((\Z^n,\ell_p)\) at every scale above a
fixed threshold requires a family of gaps covering the entire range.
For arbitrary \(p\), identifying such a family may require
detailed information about the distance spectrum.
When \(p\in\Z_{\ge1}\), however, the required
gaps are explicit, since
\[
\|z\|_p^p=\sum_{i=1}^n |z_i|^p\in\Z_{\ge0}
\qquad(z\in\Z^n).
\]
Consequently, for every integer \(m\ge0\), no lattice distance lies in
\[
\bigl(m^{1/p},(m+1)^{1/p}\bigr),
\]
whether or not either endpoint is realized. Combining the density estimate in
\cref{lem:coordinate-lattice-density} with the Jung estimate in
\cref{lem:lp-jung}, the proof below gives an explicit index \(M_{p,n}\) such
that the hypotheses of \cref{prop:distance-gap-contractibility} hold for every
such gap with \(m\ge M_{p,n}\). The resulting intervals of open scales cover
every scale above \(M_{p,n}^{1/p}\).

\begin{corollary}[Distance-gap refinement for integer
exponents]
\label{cor:integer-p-refinement}
Let \(n\ge1\) and \(p\in\Z_{\ge1}\). With \(j_{p,n}\) and \(r_{p,n}\) as in
\cref{thm:lp-lattice}, set
\[
M_{p,n}:=
\max\left\{
n,\,
\left\lfloor
r_{p,n}^{\,p}-\frac1{1-j_{p,n}}
\right\rfloor+1
\right\}.
\]
Then \(\VRlt{\Z^n,\ell_p}{r}\) is contractible whenever
\(r>M_{p,n}^{1/p}\), and \(\VRle{\Z^n,\ell_p}{r}\) is contractible
whenever \(r\ge M_{p,n}^{1/p}\).
If \(n\ge2\), then \(M_{p,n}^{1/p}<r_{p,n}\).
\end{corollary}

\begin{proof}
Write \(j:=j_{p,n}\), and let
\(\delta_p:=n^{1/p}/2\), so that \(\Z^n\) is \(\delta_p\)-dense in
\(\ell_p^n\). Set \(r_0:=r_{p,n}\). By the definition of \(r_{p,n}\),
\[
\delta_p=(1-j)r_0.
\]
We first show that every integer \(m\ge M_{p,n}\)
determines a gap to which
\cref{prop:distance-gap-contractibility} applies. Set
\[
a:=m^{1/p},
\qquad
b:=(m+1)^{1/p}.
\]
The interval \((a,b)\) contains no lattice distance, since the \(p\)-th
power of such a distance would be an integer strictly between \(m\) and
\(m+1\), which is impossible. Moreover, \(m\ge M_{p,n}\ge n\), and hence
\[
b=(m+1)^{1/p}>n^{1/p}=2\delta_p.
\]

We next prove
\[
ja+\delta_p<b.
\tag{*}
\]
We consider two cases. If \(r_0\le b\), then it follows from
\(a<b\) and \(0<j<1\) that
\[
ja+\delta_p
=ja+(1-j)r_0
<b.
\]
Now suppose that \(r_0>b\), equivalently \(r_0^p>m+1\). Since
\(m<m+1<r_0^p\), the concavity of \(x\mapsto x^{1/p}\) gives
\[
\frac{r_0-a}{r_0^p-m}\le b-a.
\]
Moreover, \(m\ge M_{p,n}\) implies
\(
m>r_0^p-\frac{1}{1-j},
\)
and hence
\[
(1-j)(r_0^p-m)<1.
\]
Combining these inequalities yields
\[
(1-j)(r_0-a)
\le (1-j)(r_0^p-m)(b-a)
<b-a.
\]
Therefore,
\[
ja+\delta_p
=ja+(1-j)r_0
=a+(1-j)(r_0-a)
<b.
\]
This proves \((*)\).

Since \(J_{\mathrm{fin}}(\ell_p^n)\le j\), inequality \((*)\) gives
\(
J_{\mathrm{fin}}(\ell_p^n)a+\delta_p
\le ja+\delta_p
<b.
\)
We have also shown that \(\Z^n\) is \(\delta_p\)-dense, that no lattice
distance lies in \((a,b)\), and that \(2\delta_p<b\). Thus all the hypotheses
of \cref{prop:distance-gap-contractibility} are satisfied with
\(X=\ell_p^n\) and \(Y=\Z^n\). Therefore
\[
\VRle{\Z^n,\ell_p}{a}
=
\VRlt{\Z^n,\ell_p}{b}
\simeq *
\]
for every integer \(m\ge M_{p,n}\).

These gapwise conclusions cover every open scale above
\(M_{p,n}^{1/p}\). Indeed, fix \(r>M_{p,n}^{1/p}\), and set
\[
m:=\lceil r^p\rceil-1.
\]
Then \(m\ge M_{p,n}\) and \(m^{1/p}<r\le(m+1)^{1/p}\). The interval conclusion in
\cref{prop:distance-gap-contractibility} therefore shows that
\(\VRlt{\Z^n,\ell_p}{r}\) is contractible. This proves the open assertion,
and \cref{prop:open-closed-threshold-equivalence} gives the closed assertion.

Assume that \(n\ge2\). Both terms defining \(j\) are greater than \(1/2\).
For \(n^{1/p}/2\), this is immediate. If \(p=1\), the first term is
\(n/(n+1)>1/2\). If \(p\ge2\), then \(s_p=p\), and
\[
\left[
2^{-1/p}\left(\frac{n}{n+1}\right)^{1/p}
\right]^p
=\frac{n}{2(n+1)}
>\frac1{2^p},
\]
because \(2^{p-1}n>n+1\). Hence \(j>1/2\), which gives
\(r_0>n^{1/p}\), so \(r_0^p>n\), and \(1/(1-j)>1\). Moreover,
\[
\left\lfloor r_0^p-\frac1{1-j}\right\rfloor+1<r_0^p.
\]
Both terms defining \(M_{p,n}\) are therefore smaller than \(r_0^p\), which
proves \(M_{p,n}^{1/p}<r_{p,n}\).
\end{proof}

Specializing \cref{cor:integer-p-refinement} to \(p=1\) makes
\(M_{p,n}\) explicit and gives the following \(\ell_1\)-lattice bound.

\integerlatticecontractibility*

\begin{proof}%
For \(p=1\), one has \(p'=\infty\), so
\cref{eq:j_{p,n}} gives
\[
j_{1,n}
=
\min\left\{\frac{n}{n+1},\frac n2\right\}
=
\frac{n}{n+1}.
\]
Thus
\[
r_{1,n}=\frac{n(n+1)}2,
\qquad
\frac1{1-j_{1,n}}=n+1
\qquad\text{and}\qquad
r_{1,n}-\frac1{1-j_{1,n}}
=\frac{n(n-1)}2-1.
\]
Therefore \(M_{1,n}=K_n\), and the result follows from
\cref{cor:integer-p-refinement}.
\end{proof}

\begin{remark}[Comparison with previous lattice bounds]
For the \(\ell_1\)-metric, Zaremsky~\cite{zaremsky} proved that
\(\VRle{\Z^n,\ell_1}{r}\) is contractible for \(r\ge n^2+n-1\). The
sufficient scale \(K_n\) in
\cref{thm:integer-lattice} satisfies
\[
\frac{K_n}{n^2+n-1}\longrightarrow\frac12
\qquad\text{as }n\to\infty.
\]

For the Euclidean metric, \cref{thm:lp-lattice} gives
\[
    j_{2,n}=\sqrt{\frac{n}{2(n+1)}},
    \qquad
    r_{2,n}
    =\frac{\sqrt n}{2\left(1-\sqrt{n/(2(n+1))}\right)}.
\]
McCarty's earlier sufficient bound for \(\VRle{\Z^n,\ell_2}{r}\) was
\(r\ge(2+\sqrt2)\sqrt n\)~\cite[Theorem~3.2]{mccarty}. The
integer-exponent refinement replaces \(r_{2,n}\) by \(M_{2,n}^{1/2}\),
which is strictly smaller when \(n\ge2\), and
\[
M_{2,n}^{1/2}
\le r_{2,n}
<\frac{2+\sqrt2}{2}\sqrt n
\]
for every \(n\ge1\). Moreover,
\[
\frac{r_{2,n}^2}{n}\longrightarrow
\left(\frac{2+\sqrt2}{2}\right)^2>1,
\qquad
\frac{1}{1-j_{2,n}}\longrightarrow2+\sqrt2,
\]
so the second term in the definition of \(M_{2,n}\) eventually dominates the
first, and \(M_{2,n}/r_{2,n}^2\longrightarrow1\). Hence
\[
\frac{M_{2,n}^{1/2}}{(2+\sqrt2)\sqrt n}
\longrightarrow\frac12
\qquad\text{as }n\to\infty.
\]
\end{remark}

\subsection{\texorpdfstring{Exact thresholds in dimensions
            one and two}{Exact thresholds in dimensions one and two}}
\label{subsec:exact-lattice-thresholds}

We next determine the exact open and closed thresholds for
\((\Z^n,\ell_p)\) in dimensions one and two. In dimension one, all the
\(\ell_p\)-metrics agree. In dimension two, the cases \(p=1\) and
\(p=\infty\) follow from the exact results already proved for those metrics.
For \(1<p<\infty\), \cref{thm:normed-lattice} first gives contractibility
for \(r\ge\frac32\,2^{1/p}\). To reach the claimed threshold
\(2^{1/p}\), we use two short-vector comparisons. At scale \(2^{1/p}\),
the closed \(\ell_p\)-Rips complex agrees with the closed
\(\ell_\infty\)-Rips complex at scale \(1\); at scale \(2\), it agrees with
the closed \(\ell_1\)-Rips complex at scale \(2\). Distance gaps extend these
two contractible complexes until they meet the range covered by
\cref{thm:normed-lattice}.

\begin{corollary}[The integer lattice in dimensions one and two]
\label{cor:lp2-exact-threshold}
Let \(n\in\{1,2\}\) and \(1\le p\le\infty\), with the convention
\(n^{1/\infty}:=1\). Then \(\VRle{\Z^n,\ell_p}{r}\) is contractible if and
only if
\[
r\ge n^{1/p}.
\]
Similarly, \(\VRlt{\Z^n,\ell_p}{r}\) is contractible if and only if
\(r>n^{1/p}\).
\end{corollary}

\begin{proof}
For \(n=1\), all the \(\ell_p\)-metrics agree with the usual metric on
\(\Z\). For finite \(p\), one has \(j_{p,1}=1/2\), so
\cref{thm:lp-lattice} gives contractibility of the open complex for \(r>1\)
and of the closed complex for \(r\ge1\); the same conclusions hold for
\(p=\infty\). The closed complex is discrete for \(r<1\), and thus not contractible; similar reasoning applies to the open complex. This proves the result for \(n=1\).

Assume now that \(n=2\). It suffices to prove the closed assertion, since the open assertion then follows from
\cref{prop:open-closed-threshold-equivalence}.
For \(p=1\),
\cref{thm:integer-lattice} gives contractibility of the closed complex for
\(r\ge2\). For
\(p=\infty\), the exact statement in \cref{thm:lp-lattice} gives
contractibility of the closed complex for \(r\ge1\).
It remains to consider \(1<p<\infty\).
We prove that
\(\VRle{\Z^2,\ell_p}{r}\) is contractible for every \(r\ge2^{1/p}\) by
considering three cases.

\emph{Case 1: \(r\ge\frac32\,2^{1/p}\).}
Coordinatewise rounding shows that
\(\Z^2\) is \(2^{1/p}/2\)-dense in \(\ell_p^2\), while
\cref{lem:bohnenblust-general} gives
\(J_{\mathrm{fin}}(\ell_p^2)\le2/3\). Thus
\cref{thm:normed-lattice} gives contractibility of the closed complex whenever
\(
r\ge\frac32\,2^{1/p}.
\)

\emph{Case 2: \(2^{1/p}\le r<2\).}
For \(z=(z_1,z_2)\in\Z^2\), integrality and \(p>1\) give
\[
\lVert z\rVert_p\le2^{1/p}
\iff |z_1|^p+|z_2|^p\le2
\iff |z_1|,|z_2|\le1
\iff \lVert z\rVert_\infty\le1.
\]
Moreover, if \(\lVert z\rVert_p>2^{1/p}\), then
\(\lVert z\rVert_\infty>1\), so integrality gives
\(\lVert z\rVert_\infty\ge2\) and hence \(\lVert z\rVert_p\ge2\).
Thus no lattice distance lies in \((2^{1/p},2)\). Consequently, for
\(2^{1/p}\le r<2\),
\[
\VRle{\Z^2,\ell_p}{r}
=
\VRle{\Z^2,\ell_p}{2^{1/p}}
=
\VRle{\Z^2,\ell_\infty}{1}
\simeq *,
\]
where the final equivalence follows from \cref{thm:lp-lattice}.

\emph{Case 3: \(2\le r<(2^p+1)^{1/p}\).}
For \(z\in\Z^2\), the inequality
\(\lVert z\rVert_1\le2\) implies
\(\lVert z\rVert_p\le\lVert z\rVert_1\le2\). Now suppose that
\(\lVert z\rVert_1>2\). Since \(z\) has integer coordinates, either one
coordinate has absolute value at least \(3\), or one coordinate has absolute
value at least \(2\) and the other has absolute value at least \(1\).
Therefore
\[
\lVert z\rVert_p
\ge
\min\left\{3,(2^p+1)^{1/p}\right\}
=(2^p+1)^{1/p}>2.
\]
Here we used \(2^p+1<3^p\). It follows that
\(\lVert z\rVert_p\le2\) if and only if \(\lVert z\rVert_1\le2\).
The same estimate shows that no lattice distance lies in
\((2,(2^p+1)^{1/p})\); the upper endpoint is attained by vectors whose
absolute coordinates are \(2\) and \(1\). Consequently, for
\(2\le r<(2^p+1)^{1/p}\),
\[
\VRle{\Z^2,\ell_p}{r}
=
\VRle{\Z^2,\ell_p}{2}
=
\VRle{\Z^2,\ell_1}{2}
\simeq *,
\]
where the final equivalence follows from \cref{thm:integer-lattice}.

These three cases cover every \(r\ge2^{1/p}\), because convexity gives
\(2(3/2)^p\le 1+2^p\), and hence
\((2^p+1)^{1/p}\ge\frac32\,2^{1/p}\).
Together with the cases \(p=1\) and \(p=\infty\), this proves the closed
assertion for \(n=2\).

\emph{Necessity.}
In the closed convention, the complex is discrete for \(r<1\) and is the
square grid graph for \(1\le r<2^{1/p}\). In the open convention, it is
discrete for \(r\le1\) and is the square grid graph for
\(1<r\le2^{1/p}\). The discrete complexes are disconnected, while the square
grid graph has nonzero first homology. Hence none of these complexes is
contractible.
\end{proof}

\section{\texorpdfstring{Closed Vietoris--Rips Complexes of
Locally Finite \(\delta\)-Dense Subsets}{Closed Vietoris--Rips Complexes of
Locally Finite delta-Dense Subsets}}
\label{sec:closed-convention}

In this section, we prove \cref{thm:main-rips} for every locally finite
\(\delta\)-dense subset of a finite-dimensional normed space. For
lattices, \cref{prop:open-closed-threshold-equivalence} converts the strict open
bound into the corresponding weak closed bound because the distance spectrum is
locally finite. A general locally finite subset need not have a locally finite
distance spectrum: its pairwise distances can accumulate in a bounded
interval. We therefore compare \(\VRle{Y}{r}\) directly with a union of
closed balls in the canonical hyperconvex ambient space.

\subsection{Nag\'orko's closed-cover nerve theorem}

A family is \emph{star-countable} if each member meets at most countably many
other members, and \emph{star-finite} if each member meets only finitely many
other members. We use the following specialization of Nag\'orko's closed-cover
nerve theorem.

\begin{theorem}[Nag\'orko's closed-cover nerve theorem
{\cite[Theorem~3.3]{nagorko}}]
\label{thm:nagorko-closed-cover}
\label{lem:closed-convex-nerve}
\label{lem:closed-polyhedral-nerve}
Let \(\mathcal Q\) be a star-countable, locally finite closed cover of a normal
space \(Z\). Suppose that its nerve is locally finite dimensional and every nonempty
intersection of a subfamily of \(\mathcal Q\) is an absolute extensor for
metric spaces. Then
\[
Z\simeq\operatorname{Nerve}(\mathcal Q).
\]
\end{theorem}

This theorem applies in particular to every locally finite, star-finite family
of nonempty closed convex subsets of a normed space. Indeed, the union of such a family is metrizable and hence normal.
Star-finiteness implies
star-countability and makes the nerve of the family locally finite
dimensional.
Local finiteness ensures that a subfamily with nonempty
total intersection is finite, and every such intersection is an absolute
extensor for metric spaces by Dugundji's extension theorem
\cite[Corollary~4.2]{dugundji}. This is the form used for all closed ball
covers in this section.

\subsection{Closed ball unions in hyperconvex spaces}
{We now prove the closed analogue of \cref{prop:open-thickening}. We first
verify the finiteness properties used to apply
\cref{thm:nagorko-closed-cover}: finite dimensionality of \(X\) and local
finiteness of \(Y\) ensure that the fixed-radius closed balls in \(H\)
centered at \(Y\) form a locally finite, star-finite family. }

\begin{lemma}[Local finiteness of fixed-radius ball covers]
\label{lem:closed-ball-cover-local-finiteness}
Let \(H\) be a normed space, let \(X\subset H\) be a finite-dimensional linear
subspace, and let \(Y\subset X\) be locally finite. For every \(\rho\ge0\), the
family
\[
\{\Ble{y}{H}{\rho}\}_{y\in Y}
\]
is locally finite and star-finite. These properties are preserved after
intersecting every member with a fixed subset of \(H\).
\end{lemma}

\begin{proof}
Fix \(z\in H\). If the neighborhood \(\Blt{z}{H}{1}\) meets
\(\Ble{y}{H}{\rho}\), then the triangle inequality gives
\(\|y\|<\|z\|+1+\rho\). Hence all such centers $y$ lie in the compact ball
\(\Ble{0}{X}{\|z\|+1+\rho}\), which contains only finitely many points of
\(Y\). This proves local finiteness.

For fixed \(y\in Y\), an intersection
\(\Ble{y}{H}{\rho}\cap\Ble{y'}{H}{\rho}\neq\emptyset\) implies
\(\|y-y'\|\le2\rho\). The compact ball
\(\Ble{y}{X}{2\rho}\) contains only finitely many points of \(Y\), which
proves star-finiteness.

Intersecting the family with a fixed subset can neither cause a neighborhood
to meet a new member nor create a new pairwise intersection. Thus both
properties are preserved.
\end{proof}

\begin{proposition}[Closed Vietoris--Rips complexes as ball unions]
\label{prop:nerve-identification}
Let \(H\) be a hyperconvex normed space, let \(X\subset H\) be a
finite-dimensional linear subspace, let \(Y\subset X\) be nonempty and locally
finite, and let \(\rho\ge0\). Then \(\VRle{Y}{2\rho}\) is the nerve of the
closed ball family \(\{\Ble{y}{H}{\rho}\}_{y\in Y}\), and
\[
\VRle{Y}{2\rho}\simeq
\bigcup_{y\in Y}\Ble{y}{H}{\rho}.
\]
\end{proposition}

\begin{proof}
If a finite nonempty set \(W\subset Y\) indexes a nonempty ball intersection,
then the triangle inequality gives \(\operatorname{diam}(W)\le2\rho\).
Conversely, if \(\operatorname{diam}(W)\le2\rho\), the balls
\(\{\Ble{w}{H}{\rho}\}_{w\in W}\) satisfy the pairwise hyperconvexity
inequalities and therefore have nonempty total intersection. This proves the
nerve identity.

By \cref{lem:closed-ball-cover-local-finiteness}, the cover  $\{\Ble{y}{H}{\rho}\}_{y\in Y}$ is
locally finite and star-finite. Its members and their intersections are closed
and convex, so \cref{thm:nagorko-closed-cover} identifies the ball union with
its nerve.
\end{proof}

\subsection{Nerve-preserving restriction to a closed tube}

The closed ball union in \cref{prop:nerve-identification} need not be convex,
so we restrict the cover to a closed tube about \(X\). The
nerve comparison is the closed-ball analogue of
\cref{lem:ball-nerve-restriction}.
Finite dimensionality ensures that
\(R_X(W)\) is attained, so the condition
\(R_X(W)\le\rho+\beta\) provides a center even when equality holds. Since the
balls are closed, hyperconvexity can be applied at the original radii
\(\rho\) and \(\beta\); the \(\varepsilon\)-shrinking used for open balls is
therefore unnecessary.

\begin{proposition}[Nerve-preserving closed-tube restriction under a finite Jung bound]
\label{prop:closed-tube-comparison}
Let \(H\) be a hyperconvex normed space, let \(X\subset H\) be a
finite-dimensional linear subspace, and let \(Y\subset X\) be nonempty and
locally finite. Fix \(\rho>0\) and \(\beta\ge0\), and suppose that
\[
J_{\mathrm{fin}}(X)\le\frac{\rho+\beta}{2\rho}.
\]
Then the covers
\[
\{\Ble{y}{H}{\rho}\}_{y\in Y}
\quad\text{and}\quad
\{\Ble{y}{H}{\rho}\cap\Tle{X}{H}{\beta}\}_{y\in Y}
\]
have the same nerve. Consequently,
\[
\bigcup_{y\in Y}\Ble{y}{H}{\rho}
\simeq
\left(\bigcup_{y\in Y}\Ble{y}{H}{\rho}\right)
\cap\Tle{X}{H}{\beta}.
\]
\end{proposition}

\begin{proof}
Restriction to the closed tube gives one inclusion of nerves. For the reverse
inclusion, let \(W\subset Y\) index a nonempty intersection of the closed
\(\rho\)-balls. Then \(\operatorname{diam}(W)\le2\rho\). If
\(\operatorname{diam}(W)=0\), then \(W=\{w\}\), and we take \(c=w\). If
\(\operatorname{diam}(W)>0\), finite dimensionality gives a circumcenter
\(c\in X\). Since \(X\) carries the metric induced from \(H\),
\[
\max_{w\in W}d_H(c,w)
=R_X(W)
\le J_{\mathrm{fin}}(X)\cdot\operatorname{diam}(W)
\le \rho+\beta.
\]
Thus, in either case,
\[
    d_H(w,w')\le2\rho,
    \qquad
    d_H(c,w)\le\rho+\beta
\]
for all \(w,w'\in W\). The closed balls
\(\{\Ble{w}{H}{\rho}\}_{w\in W}\), together with
\(\Ble{c}{H}{\beta}\), therefore satisfy the pairwise hyperconvexity
inequalities. Hyperconvexity gives a point common to all these balls. Since
\(c\in X\), this point lies in \(\Tle{X}{H}{\beta}\), so the restricted
balls indexed by \(W\) also have nonempty intersection. The two covers have
the same nerve.

The set \(\Tle{X}{H}{\beta}\) is closed and convex. By
\cref{lem:closed-ball-cover-local-finiteness}, both covers are locally finite
and star-finite, and all their nonempty intersections are closed and convex.
Applying \cref{thm:nagorko-closed-cover} to both covers gives the stated
homotopy equivalence.
\end{proof}

\subsection{Density fills the closed tube}
{
We now combine the closed tube restriction with the \(\delta\)-density of \(Y\)
to prove \cref{thm:main-rips}. 
}

\closedellonecontractibility*

\begin{proof}%
If \(X=\{0\}\), then \(\delta\)-density forces \(Y=\{0\}\), and the result is
immediate. Assume \(\dim X=n\ge1\), and put
\[
J:=J_{\mathrm{fin}}(X).
\]
Then \(1/2\le J<1\): the lower bound follows from two-point sets, and the upper
bound follows from \cref{lem:bohnenblust-general}. Set
\[
\rho:=\frac r2,
\qquad
\beta:=(2J-1)\rho.
\]
Then \(\beta\ge0\),
\[
J=\frac{\rho+\beta}{2\rho},
\qquad
\delta+\beta\le\rho.
\]

Let \(H:=\ell_\infty(B_{X^*})\), and let
\(\iota:X\to H\) be the canonical isometric embedding. Write
\[
T:=\Tle{\iota(X)}{H}{\beta}.
\]
Because \(\iota(X)\) is finite dimensional, the distance from any point of
\(H\) to \(\iota(X)\) is attained: a minimizing sequence is bounded and hence
has a convergent subsequence in \(\iota(X)\). Consequently,
\[
T=\iota(X)+\Ble{0}{H}{\beta}.
\]
Since \(\iota(X)\) and \(\Ble{0}{H}{\beta}\) are convex, the tube \(T\) is
convex.

The space \(H\) is hyperconvex by
\cite[Example~2.22 and Proposition~2.18]{limmemoliokutan}. Since \(\iota\) is
an isometry,
\(J_{\mathrm{fin}}(\iota(X))=J_{\mathrm{fin}}(X)=J\), and
\(\iota(Y)\) is locally finite in \(\iota(X)\). Using \(\iota\) to identify
\(\VRle{Y}{2\rho}\) with \(\VRle{\iota(Y)}{2\rho}\),
\cref{prop:nerve-identification,prop:closed-tube-comparison} give
\[
\VRle{Y}{2\rho}
\simeq \bigcup_{y\in Y}\Ble{\iota(y)}{H}{\rho}
\simeq
\left(\bigcup_{y\in Y}\Ble{\iota(y)}{H}{\rho}\right)\cap T.
\]

It remains to identify the restricted union with \(T\). Since the restricted
union is contained in \(T\) by definition, it suffices to prove
\[
T\subset\bigcup_{y\in Y}\Ble{\iota(y)}{H}{\rho}.
\]
Let \(z\in T\). Choose \(x\in X\) with
\(d_H(z,\iota(x))\le\beta\), and then choose \(y\in Y\) with
\(d_X(x,y)\le\delta\).
Since \(\iota\) is an isometry,
\[
\begin{aligned}
    d_H(z,\iota(y))
     & \le d_H(z,\iota(x))+d_H(\iota(x),\iota(y)) \\
     & =d_H(z,\iota(x))+d_X(x,y)                  \\
     & \le\beta+\delta\le\rho.
\end{aligned}
\]
Thus \(z\in\Ble{\iota(y)}{H}{\rho}\), and hence
\[
\left(\bigcup_{y\in Y}\Ble{\iota(y)}{H}{\rho}\right)\cap T=T.
\]
The space \(T\) is convex and therefore contractible. Since \(r=2\rho\), the
result follows.
\end{proof}

Bohnenblust's estimate
\(J_{\mathrm{fin}}(X)\le n/(n+1)\) gives the following
dimension-dependent consequence of \cref{thm:main-rips}.

\begin{corollary}[{Finite-dimensional bound for closed
Vietoris--Rips complexes}]
\label{cor:finite-dimensional-closed}
Let \(X\) be an \(n\)-dimensional normed space with \(n\ge1\), and let
\(Y\subset X\) be locally finite and \(\delta\)-dense. Then \(\VRle{Y}{r}\) is
contractible for every
\[
r\ge(n+1)\delta.
\]
\end{corollary}

\begin{proof}
By \cref{lem:bohnenblust-general},
\(J_{\mathrm{fin}}(X)\le n/(n+1)\). Hence
\(\delta/(1-J_{\mathrm{fin}}(X))\le(n+1)\delta\), and
\cref{thm:main-rips} applies.
\end{proof}

For \(1\le p\le\infty\), the exact thresholds in
\cref{cor:lp2-exact-threshold} concern the lattice \(\Z^2\).
For arbitrary locally finite \(\delta\)-dense subsets of
\((\R^2,\ell_p)\), \cref{thm:main-rips} gives a closed contractibility bound
in terms of \(J_{\mathrm{fin}}(\ell_p^2)\). To make this bound explicit, set
\begin{equation}
\label{eq:planar-lp-constant}
c_p:=
\begin{cases}
2,                     & p\in\{1,\infty\}, \\[2mm]
\dfrac{1}{1-j_{p,2}}, & 1<p<\infty.
\end{cases}
\end{equation}
By \cref{lem:lp-jung,lem:herrlich-embedding},
\begin{equation}
\label{eq:planar-jung-bound}
J_{\mathrm{fin}}(\ell_p^2)\le1-\frac1{c_p}
\qquad(1\le p\le\infty).
\end{equation}
In particular, \(c_2=(3+\sqrt3)/2\), improving the coefficient \(3\) given
by \cref{cor:finite-dimensional-closed} in the Euclidean plane.

\begin{corollary}[{Closed bounds in \(\ell_p^2\)}]
\label{prop:two-dimensional-offset-comparison}
\label{cor:closed-planar-lp}
Let \(1\le p\le\infty\), let \(c_p\) be defined by
\eqref{eq:planar-lp-constant}, and let
\(Y\subset(\R^2,\ell_p)\) be locally finite and \(\delta\)-dense. Then
\(\VRle{Y}{r}\) is contractible for every \(r\ge c_p\delta\).
\end{corollary}

\begin{proof}
By \eqref{eq:planar-jung-bound},
\[
\frac{\delta}{1-J_{\mathrm{fin}}(\ell_p^2)}\le c_p\delta.
\]
The assertion follows from \cref{thm:main-rips}.
\end{proof}

We conclude by discussing the sharpness of the coefficient \(c_p\) in
\cref{prop:two-dimensional-offset-comparison}. For \(1\le p\le\infty\),
\cref{lem:coordinate-lattice-density} shows that \(\Z^2\) is
\(\delta_p\)-dense in \(\ell_p^2\), where
\(\delta_p:=2^{1/p}/2\) and \(2^{1/\infty}:=1\). By
\cref{cor:lp2-exact-threshold}, its exact threshold is
\(2^{1/p}=2\delta_p\). Rescaling the lattice therefore
shows that no coefficient below \(2\) can hold uniformly for all locally
finite \(\delta\)-dense subsets of \((\R^2,\ell_p)\).
Since \(c_1=c_\infty=2\), the
planar bounds are sharp for \(p=1\) and \(p=\infty\); see
\cref{cor:lp-open-dense,cor:open-l12-sharp,cor:closed-planar-lp}.
The
coefficient \(2\) is also optimal in dimension one. We do not claim that
\(c_p\) is optimal for \(1<p<\infty\).

\section*{Generative AI Disclosure}

The authors developed the argument by restricting the cover to a convex tube
around the embedded normed space. During the preparation of the manuscript,
they used ChatGPT, Claude Code, and Codex to discuss the scope and quantitative
limits of the argument in the \(\ell_p\)-lattice setting, including the
distance-gap refinement in \cref{cor:integer-p-refinement}, and to help
identify relevant literature, including McCarty's
dissertation~\cite{mccarty}. The authors take full responsibility for the
contents of the paper.

\bibliographystyle{amsplain}
\bibliography{arxiv_2}

\end{document}